\documentclass[12pt,twoside]{article}
\usepackage{graphicx}
\usepackage{tikz}
\usepackage{algorithm}
\usepackage{algorithmic}

\newtheorem{theorem}{Theorem}
\newtheorem{definition}[theorem]{Definition}
\newtheorem{lemma}[theorem]{Lemma}
\newtheorem{corollary}[theorem]{Corollary}
\newenvironment{keywords}{\centerline{\bf\small
Keywords}\begin{quote}\small}{\par\end{quote}\vskip 1ex}
\def\expec#1#2{{\rm I\!E}_{#1}\left[#2\right]}
\def\vari#1#2{\mbox{Var}_{#1}\left[#2\right]}
\def\covi#1#2#3{\mbox{Cov}_{#1}\left[#2,\,#3\right]}
\def\nq{\hspace{-1em}}
\def\SetN{I\!\!N}
\def\SetR{I\!\!R}
\def\qmbox#1{{\quad\mbox{#1}\quad}}

\def\req#1{(\ref{#1})}
\def\qed{}
\def\beq{\begin{equation}}    \def\eeq{\end{equation}}
\def\beqn{\begin{displaymath}}\def\eeqn{\end{displaymath}}
\def\bqa{\begin{eqnarray}}    \def\eqa{\end{eqnarray}}
\def\bqan{\begin{eqnarray*}}  \def\eqan{\end{eqnarray*}}
\def\a{\alpha}
\def\b{\beta}
\def\d{\delta}

\begin{document}


\title{\vspace{-4ex}
\vskip 2mm\bf\Large\hrule height5pt \vskip 4mm
A Bayesian View of the \\ Poisson-Dirichlet Process
\vskip 4mm \hrule height2pt}

\author{{\bf Wray~Buntine}\\[1mm]
{\tt wray.buntine@nicta.com.au}\\
NICTA and Australian National University\\
Locked Bag 8001, Canberra ACT 2601, Australia
\and
{\bf Marcus~Hutter}\\[1mm]
{\tt marcus.hutter@anu.edu.au}\\
Australian National University and NICTA\\
RSISE, Daley Road, Canberra ACT 0200, Australia
}

\date{15 February 2012}
\maketitle

\newpage

\begin{abstract}
The two parameter Poisson-Dirichlet Process (PDP),
a generalisation of the Dirichlet Process, is increasingly
being used for probabilistic modelling in discrete areas such as
language technology, bioinformatics, and image analysis.
There is a rich literature about the PDP and its
derivative distributions such as the Chinese Restaurant Process.
This article reviews some of the basic theory
and then the major results needed for Bayesian
modelling of discrete problems
including details of priors, posteriors and computation.

The PDP is a generalisation of the Dirichlet distribution that
allows one to build distributions over partitions, both finite
and countably infinite.  The PDP has two other remarkable
properties:  first it is partially conjugate to itself,
which allows one to build hierarchies of PDPs, and
second using a marginalised relative the
Chinese Restaurant Process (CRP),
one gets fragmentation and clustering properties
that lets one layer partitions to build trees.
This article presents the basic theory for understanding
the notion of partitions and distributions over them,
the PDP and the CRP, and the important properties
of conjugacy, fragmentation and clustering, as well as some key
related properties such as consistency and convergence.
This article also presents a
Bayesian interpretation of the Poisson-Dirichlet
process: it is based on an improper and infinite dimensional
Dirichlet distribution.
This interpretation requires technicalities
of priors, posteriors and Hilbert spaces, but conceptually, this
means we can understand the process as just another Dirichlet and
thus all its sampling properties emerge naturally.

The theory of PDPs is usually presented for continuous distributions (more
generally referred to as non-atomic distributions), however, when applied
to discrete distributions its remarkable conjugacy property emerges.
This context and basic results are also presented,
as well as techniques for computing the second order Stirling
numbers that occur in the posteriors for discrete distributions.
\end{abstract}


\newpage
\tableofcontents

\newpage
\section{Introduction}

The {\em two-parameter Poisson-Dirichlet process} (PDP), also known
as the Pitman-Yor process (named so in \cite{Ish_jasa01}), is an
extension of the {\em Dirichlet process} (DP). Related is a
marginalisation known as the Chinese
Restaurant Process (CRP) which gives an elegant analogy of
incremental sampling of partitions. These models have proven useful
in a number of ways as tools for non-parametric and hierarchical
Bayesian modelling, especially in discrete domains such as with language
and images where one wants to develop hierarchical models,
or be flexible with dimension.

In language domains, PDPs are proving useful for full probability
modelling of various phenomena including n-gram modelling and
smoothing \cite{teh_acl06,Goldwater06,mochihashi_nips07}, dependency
models for grammar \cite{Johnson06,Wallach08}, and for data
compression \cite{WooArcGas2009a}. The PDP-based n-gram models
correspond well to versions of Kneser-Ney smoothing
\cite{teh_acl06}, the state of the art method in applications. These
models are intriguing from the probability perspective, as well as
sometimes being competitive in terms of performance. More
generally, the models are also being used for clustering
\cite{GreenRich01,Rasmussen2000}, and for related tasks such as
image segmentation \cite{Sudderth08},
 relational modelling \cite{XuTresp06},
and exemplar-based clustering \cite{Tarlow08}.

PDPs and their associated distributions
are basically a tool for modelling two kinds of objects:
mixture models and partitions.
They can then be used to model trees and other hierarchical structures.
Section~\ref{sct-disc} introduces the interrelated concepts
of a simple mixture model and a partition along with some of their
statistical complications.
The definition of a PDP is then presented in Section~\ref{sct-back}.
The theory is well developed for the more general context of
continuous distributions, and that theory is reviewed here.
Details of statistics, consistency and the forms of
posteriors are given in Section~\ref{sct-dd}.
Statistical analysis of partitions is given in Section~\ref{sct-pp}
which arises when one marginalises the posterior of the PDP to
obtain the CRP,
which can be further marginalised to obtain a
distribution on partition sizes.
Results on fragmentation and grouping of partitions in Section~\ref{sct-fc}
allow one to extend the CRP distribution to trees
({\it i.e.}, nested partitions).

A new Bayesian
interpretation and definition of the PDP is then given in
Section~\ref{sct-dm}. This uses the methodology of improper
priors too show that the distribution on the
infinite probability vector underlying the PDP is in fact
an infinite improper Dirichlet.
From this one can readily obtain
all the standard sampling and additivity properties for the PDP.

With the use of PDPs increasing in computer science applications,
where sophisticated discrete probabilistic modelling is required,
this article reviews and summarises the basic theory of PDPs in the
discrete context in Section~\ref{sct-post}.  This context has
distinct properties where the distribution on partition
size, introduced above, is needed to express posteriors
using Stirling numbers of the second kind.
These play the role of the Gamma function that occurs in the
posterior for a Dirichlet distribution.
Basic conjugacy results for the
PDP are then presented
that make it such an important distribution for hierarchical
Bayesian reasoning.

We have attempted where possible to follow conventions used in the
mathematical statistics community, and to provide some
pointers to that literature.

\section{Mixture Models and Partitions}\label{sct-disc}

\begin{figure*}
\begin{center}
\input{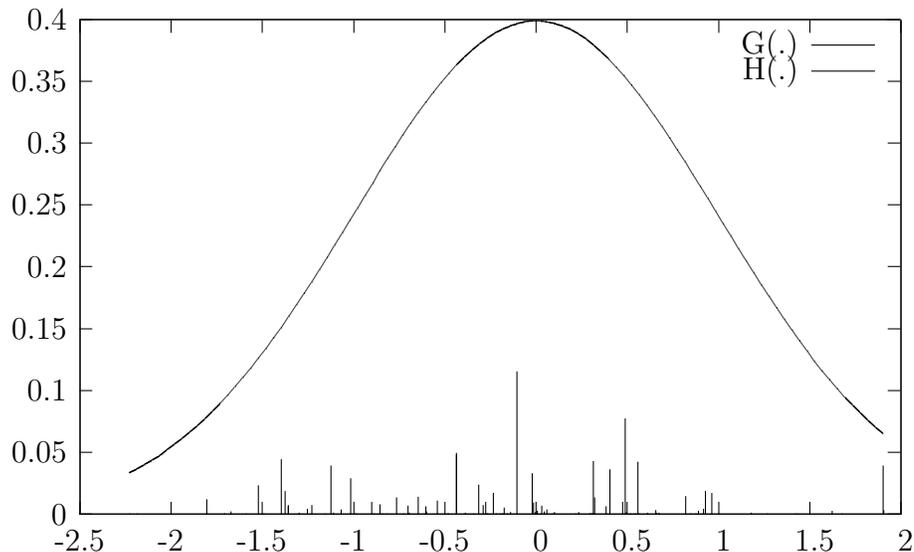}
\caption{Gaussian base distribution $H(\cdot)$ and a resultant impulse mixture model $G(\cdot)$.  The infinite number of smaller impulses do not show.}
\label{fig:ge}
\end{center}
\end{figure*}

Before introducing the PDP, we introduce the basic context of its
use, a simple mixture model.
The kinds of mixture models we consider require as input a {\it base}
probability distribution $H(\cdot)$ on a measurable space ${\cal
X}$, and yield a discrete distribution on a finite or countably
infinite subset of ${\cal X}$.  This means the output distribution
is a weighted set of impulses at points in ${\cal X}$.
\begin{definition}[Impulse mixture model]
Given a probability distribution $H(\cdot)$ on a measurable space ${\cal
X}$, assume the values
$X^*_k \in {\cal X}$ for $k=1,...,\infty$
are independently and identically distributed
according to $H(\cdot)$.
Also an infinite dimensional
probability vector $\vec{p}$ is sampled from a distribution $Q(\cdot)$
independently of each $X^*_k$
so $0\leq p_k \leq 1$ and $\sum_{k=1}^\infty p_k = 1$.
Then
\beq\label{dfn-form}
  \sum_{k=1}^\infty p_k \delta_{X^*_k}(\cdot)
\eeq
is an {\em impulse mixture model} with
{\em base distribution} $H(\cdot)$
and {\em probability distribution} $Q(\cdot)$.
Note $\delta_{X^*_k}(\cdot)$ is a
discrete measure concentrated at $X^*_k$.
\end{definition}
An example of an impulse mixture model
is given in Figure~\ref{fig:ge} where the base
distribution is a Gaussian and the probability
vector $\vec{p}$ was generated by a so-called stick-breaking method
(introduced later).

These kinds of models have a long history and appear in
various forms.
In \cite{Pitman96} the
species sampling model is developed in the
context of the PDP
and attributed to R.A.\ Fisher for simple models of
animal species.
An alternative derivation takes the
view that the $\vec{p}$ are initially
unnormalised, so-called random measures,
which one then normalises \cite{JaLiPr2009}.
There are various schemes for
sampling the probability vector $\vec{p}$ including
normalized random measures \cite{JaLiPr2009}
and general species sampling schemes \cite{IsJa2003}.
The Poisson-Dirichlet distribution we present here
is a particular version that is partially
conjugate and thus admits convenient Bayesian computation.

When the base distribution is continuous, such as with the Gaussian,
one can see that almost surely no two values $X^*_k$ and $X^*_l$
(for $k\neq l$) would be the same.
The general property is described as  {\em non-atomic}, which means
$H(X)=0$ for all $X \in \cal X$, then samples from $H(\cdot)$ are
almost surely distinct.
The counter property is
where the distribution is {\em discrete}, so $H(X)>0$ for all $X \in \cal X$
and thus it is always possible that $X^*_k=X^*_l$ for some $l \neq k$.
This distinction,
non-atomic versus discrete, has important consequences for posterior analysis
as explained in Section~\ref{sct-post}.
Mixed base distributions are not considered.

\begin{figure}
\begin{center}
\input{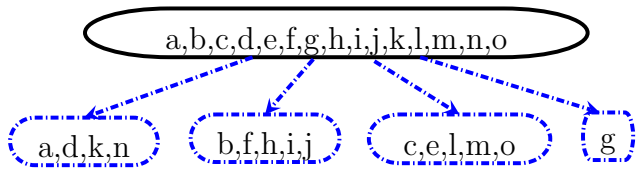}
\caption{Shows a partition of a set of letters $\{a,b,...,o\}$, illustrating a
{\em fragmentation}.  The top node is the full set and the bottom row is the partition displayed in order of least elements, {\it i.e.} $a, b, c, g$.}
\label{fig:part4}
\end{center}
\end{figure}

\begin{figure}
\begin{center}
\input{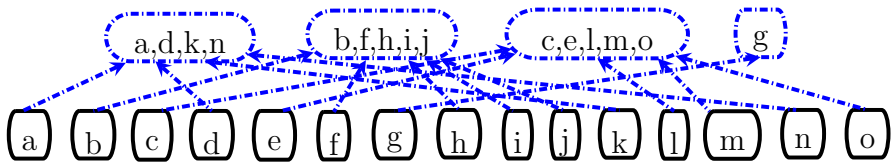}
\caption{Shows a partition applied to a set of sets of letters
 $\{\{a\},\{b\},...,\{o\}\}$, illustrating a {\em coagulation}.  The bottom set of nodes lists the sets, and the top row applies the partition to the set of sets to get subsets of sets which are then unioned.  The top row is displayed in
order of least elements.}
\label{fig:part3}
\end{center}
\end{figure}

When sampling from an impulse mixture model, one would observe a
sequence of $N$ data values
$X_1, X_2, X_3,...,X_N$,
and one may also assume corresponding indices,
the $k$ in Formula~\req{dfn-form}, also exist, so these might be
$k_1, k_2, k_3, ..., k_N$.
Since the actual indices are latent and not part of the observed
data for the model of Formula~\req{dfn-form}, the latent indices
can be arbitrarily relabelled and converted to a normal form.
Such a normal form where the labels are irrelevant is called a
{\em partition}: it defines
the grouping of items in the sequence, not the actual indices
assigned.
\begin{definition}[Partition]\label{defn-part}
A {\em partition} $P$ of a countable set $X$ is a mutually exclusive
and exhaustive set of subsets of $X$.
The {\em partition size} of $P$ is given by the number of sets $|P|$.
\end{definition}
In our case, we make a partition of the data sequence
$\{X_1, X_2, X_3,...,X_N\}$.
Partitions are important because, not only do they provide the
assignments in a mixture model, they can also be used
as a primitive in the generation of many data structures.
For instance, random trees can be created using partitions
either by a top-down process of fragmentation,
illustrated in Figure~\ref{fig:part4},
or by a bottom-up process of coagulation,
illustrated in Figure~\ref{fig:part3}.
When storing and sampling partitions, however, one may need to
order them.  Most importantly, when sampling infinite partitions
or infinite probability vectors $\vec{p}$, one needs to order
them roughly by size:  there is no point in generating a million
infinitesimal probabilities before generating one with significant size.

One obvious ordering condition to try is by size.  One could make the
probability vector  $\vec{p}$
ordered so that $p_{k+1}\leq p_k$ for all $k$,
or order bins  in the sampled partition  by their decreasing number of entries.
From a statistical perspective, however, this is impractical
because one does not know the value of each $p_{k}$, and the sizes of the
bins will vary during sampling.

The {\em size-biased order} and the {\em order of least elements}
are two ways of ordering sets
that order partitions by their order of first occurrence in
a sequence.   The sets in the partition can then be enumerated
according to this order.
For instance,
assuming the data sequence is `a', 'b', 'c', ...,
then the top
partition in Figure~\ref{fig:part3} is listed left to right
in order of least elements.
We define slightly different versions of these depending
in the representation of the the set being ordered.
\begin{definition}[Orders for sets]
Size-biased orders are defined for index sequences, probability
vectors and partitions:
\begin{itemize}
\item
An {\em index sequence} $I$ of length $N$ given by $k_1,k_n,...,k_N$
is in {\em size-biased order} if $k_1=1$ and
$k_n\leq 1+\max_{1\leq i <n}k_i$ for $n=2,...,N$.
\item
An {\em infinite probability vector} $\vec{p}$ is in {\em size-biased order}
if the elements $p_k$ have been reordered to be in their order of
first occurrence in a random sample from $\vec{p}$.
\item
A {\em partition} $P$ is in {\em order of least elements} if members are
listed in order of their least element.
\end{itemize}
\end{definition}
Note these three definitions are related as follows:
if we take an infinite  random sample from a probability vector $\vec{p}$,
and then renumber the index sequence so they are in size-biased order,
then the corresponding renumbering converts the
probability vector $\vec{p}$ to size-biased order.
If we take an infinite random sample from an impulse mixture model
with probability vector $\vec{p}$ and non-atomic base distribution
over $\cal X$, partition the sample according the values from $\cal X$,
and then label the entries by the order of occurrence,
$1,2,3,...$, then the order of least elements of the partition
also yields a size-biased order of the probability vector $\vec{p}$.

A normal form for index sequences can be derived by renumbering
indices so that the sequence is in size-biased order:
indices are renumbered to 1,2,3,... according to their
first occurrence in the sequence.
So a size-biased ordered {\em renumbering} of
the indices $12, 435, 7198, 12, 12,  35, 7198$
is $1, 2, 3, 1, 1,  4, 3$.
Also, a partition of a data sequence can be represented by
a size-biased order of the indices.  The top
partition in Figure~\ref{fig:part3},
for instance, can be represented by the
sequence ``1 2 3 1 3 2 4 2 2 2 1 3 3 1 3''.

\section{The Poisson-Dirichlet Process}\label{sct-back}

For an impulse mixture model, we
not only need a base distribution $H(\cdot)$
but also need to specify
the probability vector $\vec{p}$. Within the PDP literature,
$\vec{p}$ follows a two parameter Poisson-Dirichlet distribution
\cite{PitmanYor97}
when the probabilities $p_k$ are ordered by decreasing size,
or equivalently by the Griffiths-Engen-McCloskey (GEM) distribution
\cite{pitmanCSP}
when the probabilities $\vec{p}$ have a size-biased order.

The GEM distribution is defined via the so-called
``stick-breaking" model which goes as follows:
\begin{figure*}
\begin{center}
\input{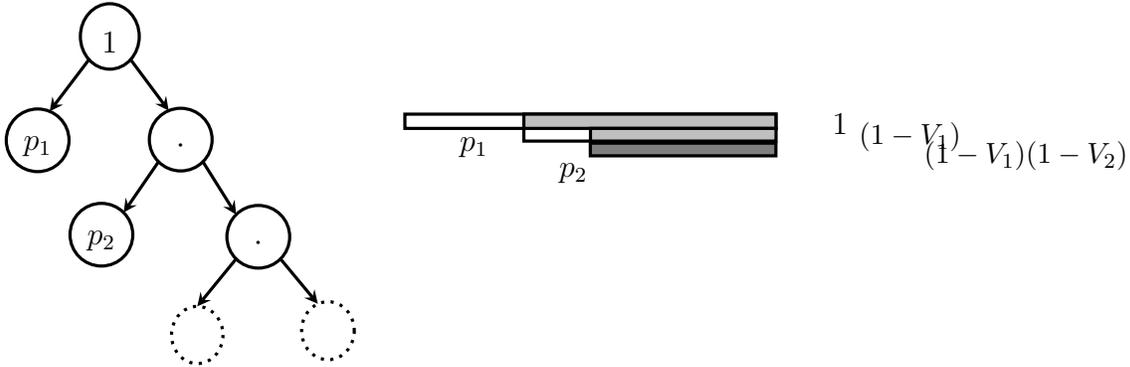}
\caption{The stick-breaking analogy for a GEM.  Left shows the tree of
probabilities being generated.  Middle shows the stick lengths, with
$p_1$ and $p_2$ already broken off, and the remainder, $(1-V_1)(1-V_2)$,
about to be broken.}
\label{fig:sticks}
\end{center}
\end{figure*}
\begin{enumerate}
\item We take a stick of length one and randomly break it into
    two parts with proportions $V_1$ and $1-V_1$. The first
    broken stick has length $V_1$.
\item We then take the remaining part, of length $1-V_1$ and
    apply the same process to randomly break into proportions
    $V_2$ and $1-V_2$. This second broken stick is the first
    part, of length $(1-V_1)V_2$.
\item Again, we take the remaining part, of length
    $(1-V_1)(1-V_2)$ and apply the same process to randomly
    partition into proportions $V_3$ and $1-V_3$. This third
    broken stick is the first part, of length
    $(1-V_1)(1-V_2)V_3$.
\item ...
\end{enumerate}
Formally, this goes as follows:
\begin{definition}[GEM distribution]\label{defn-GEM}
For $0\leq a < 1$  and $b >-a$, suppose that
independent random variables $V_k$ are such that $V_k$ has
$\mbox{Beta}( 1-a, b+ k\,a)$ distribution. Let
\beqn
  p_1 ~=~V_1,~~~~ p_k ~=~(1 - V_1) \cdots (1- V_{k-1}) V_k ~~~~k \geq 2~.
\eeqn
Define the {\em Griffiths-Engen-McCloskey distribution}
with parameters $a,b$, abbreviated $\mbox{GEM}(a,b)$ to be the
resultant distribution of $(p_1,p_2,...)$.
\end{definition}
\begin{definition}[Poisson-Dirichlet distribution]\label{defn-PY}
Let $(\tilde{p}_1,\tilde{p}_2,...)~\sim~\mbox{GEM}(a,b)$
and define $\vec{p}=(p_1,p_2,...)$ to be their sorted values
so that $p_1\geq p_2 \geq \cdots$.
Then $\vec{p}$ follows the {\em Poisson-Dirichlet distribution}
with parameters $a,b$, abbreviated $\mbox{PDD}(a,b)$.
\end{definition}
Here the parameter $a$ is usually called the {\em discount
parameter} in the literature, and $b$ is called the {\em
concentration parameter}.  The term concentration when used
in the statistics community usually means a quantity that
behaves like the inverse of a variance.

\begin{figure*}
\begin{center}
\includegraphics[width=3.in]{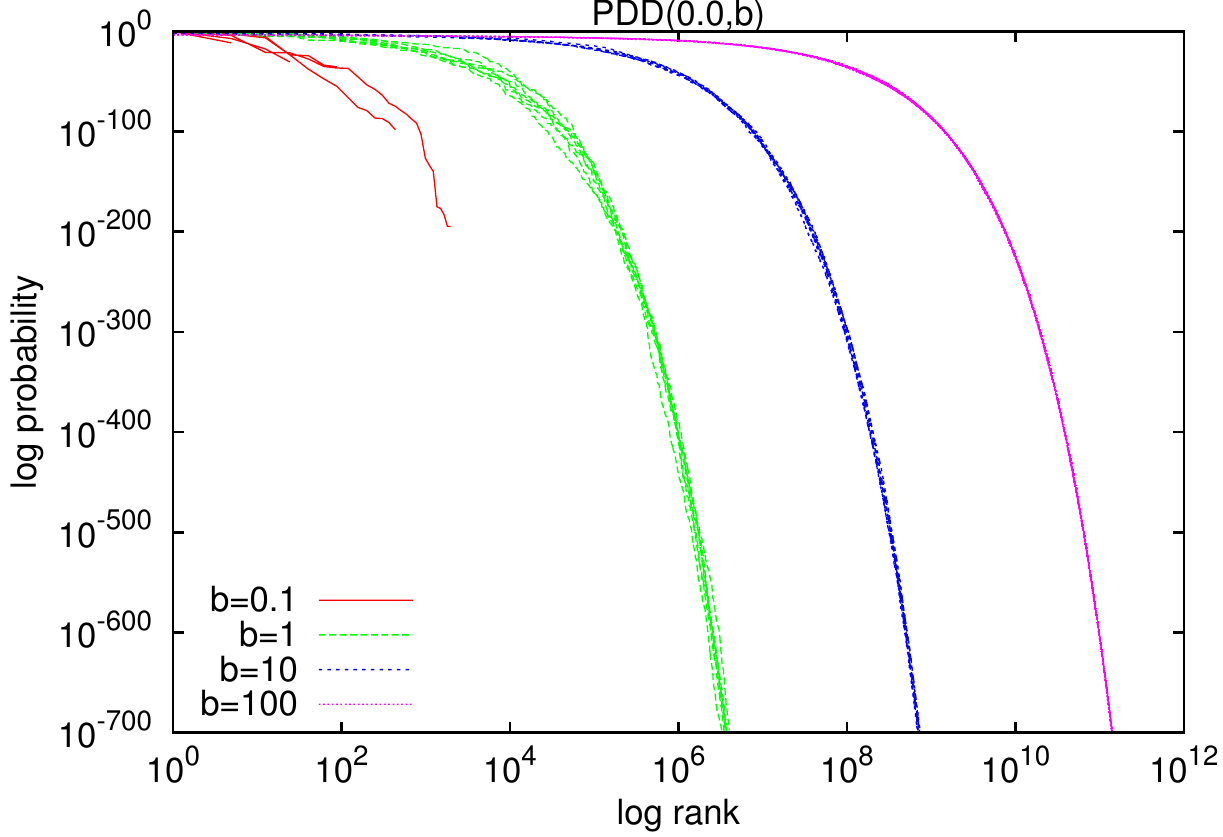}\includegraphics[width=3.in]{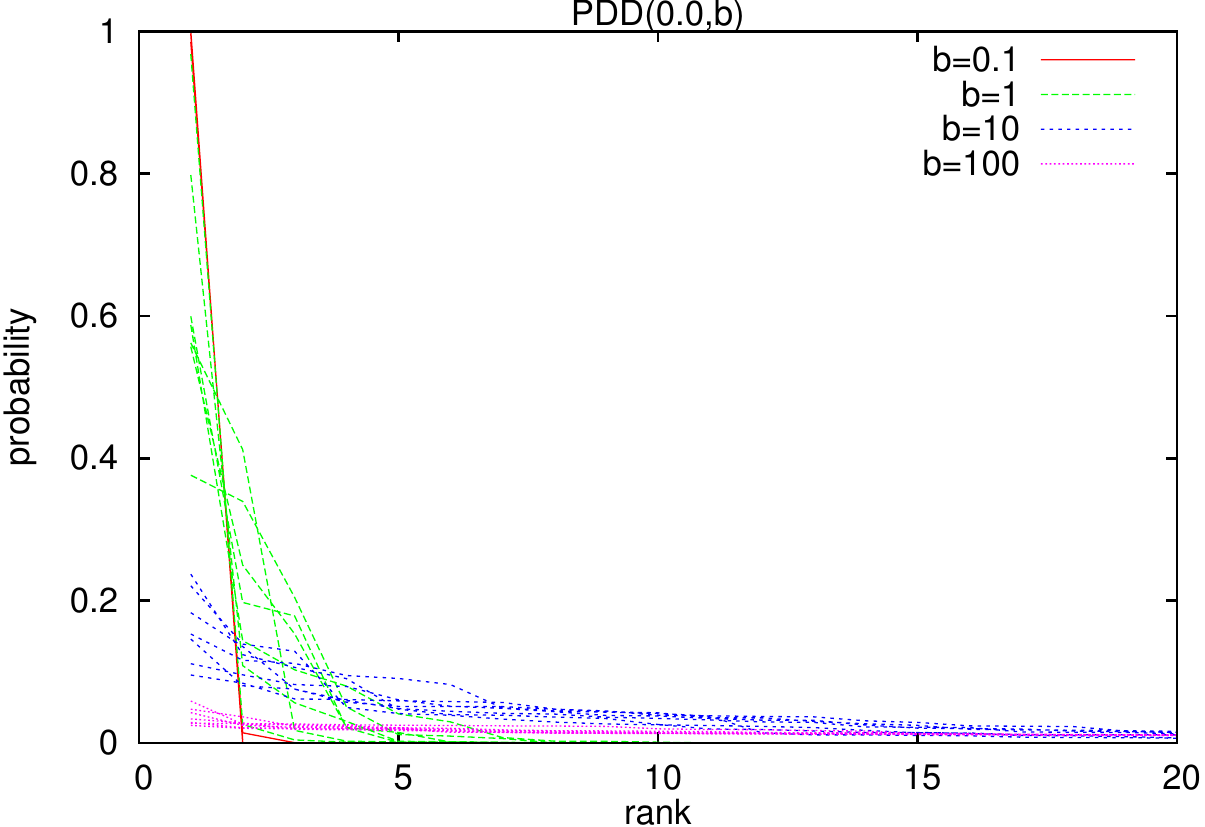}\\\vspace*{5pt}
\includegraphics[width=3.in]{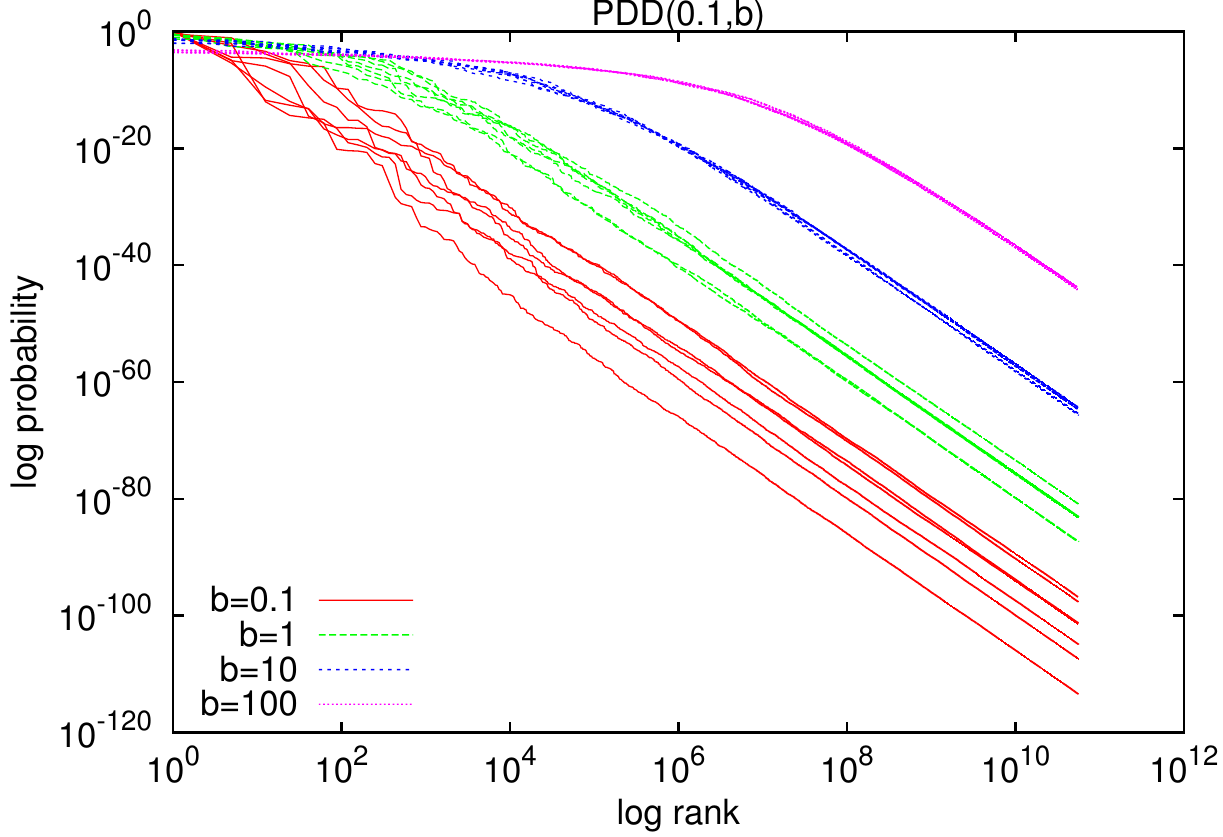}\includegraphics[width=3.in]{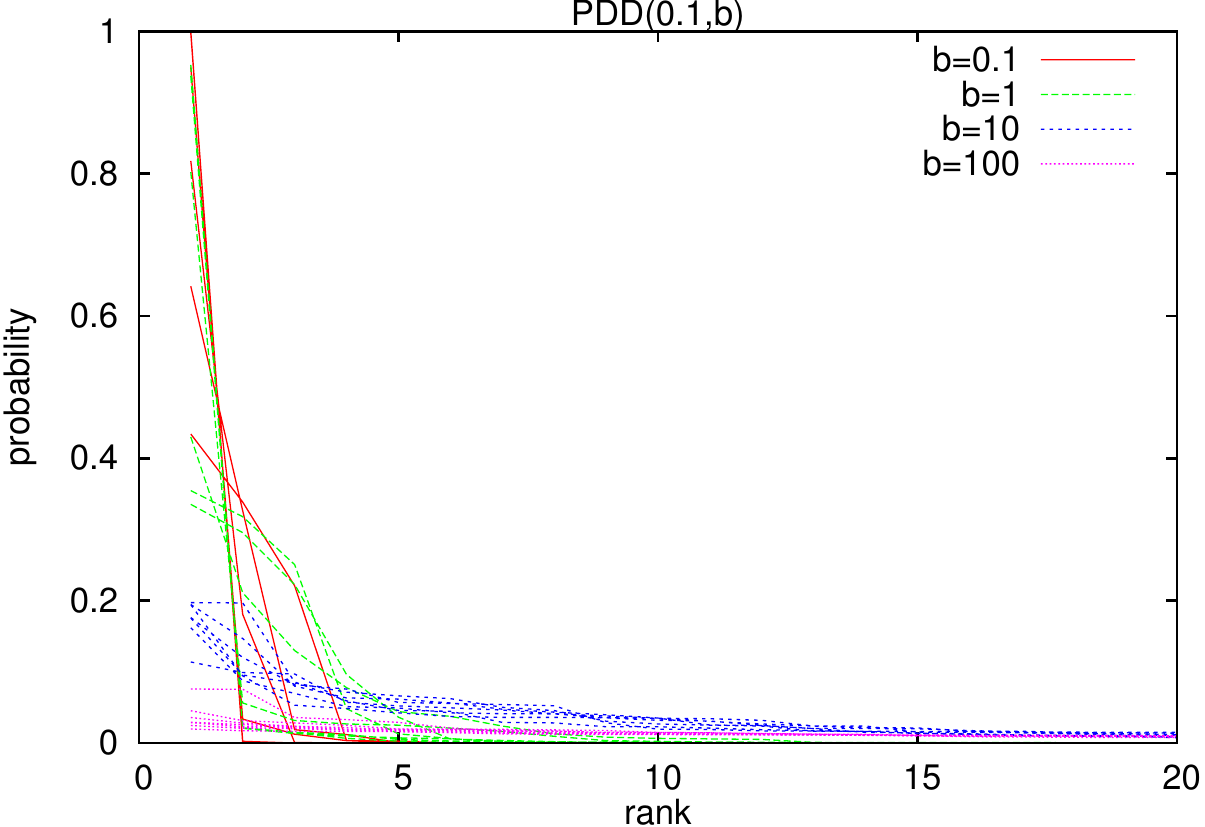}\\\vspace*{5pt}
\includegraphics[width=3.in]{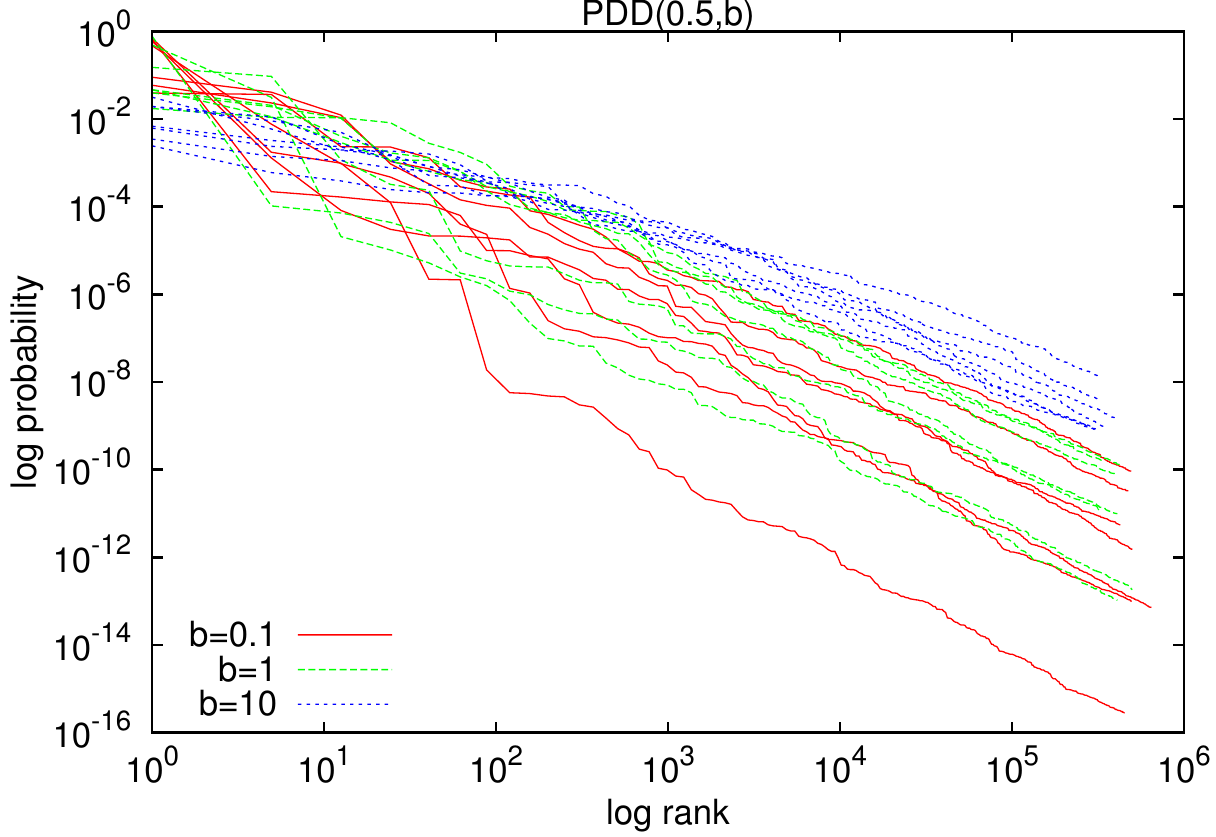}\includegraphics[width=3.in]{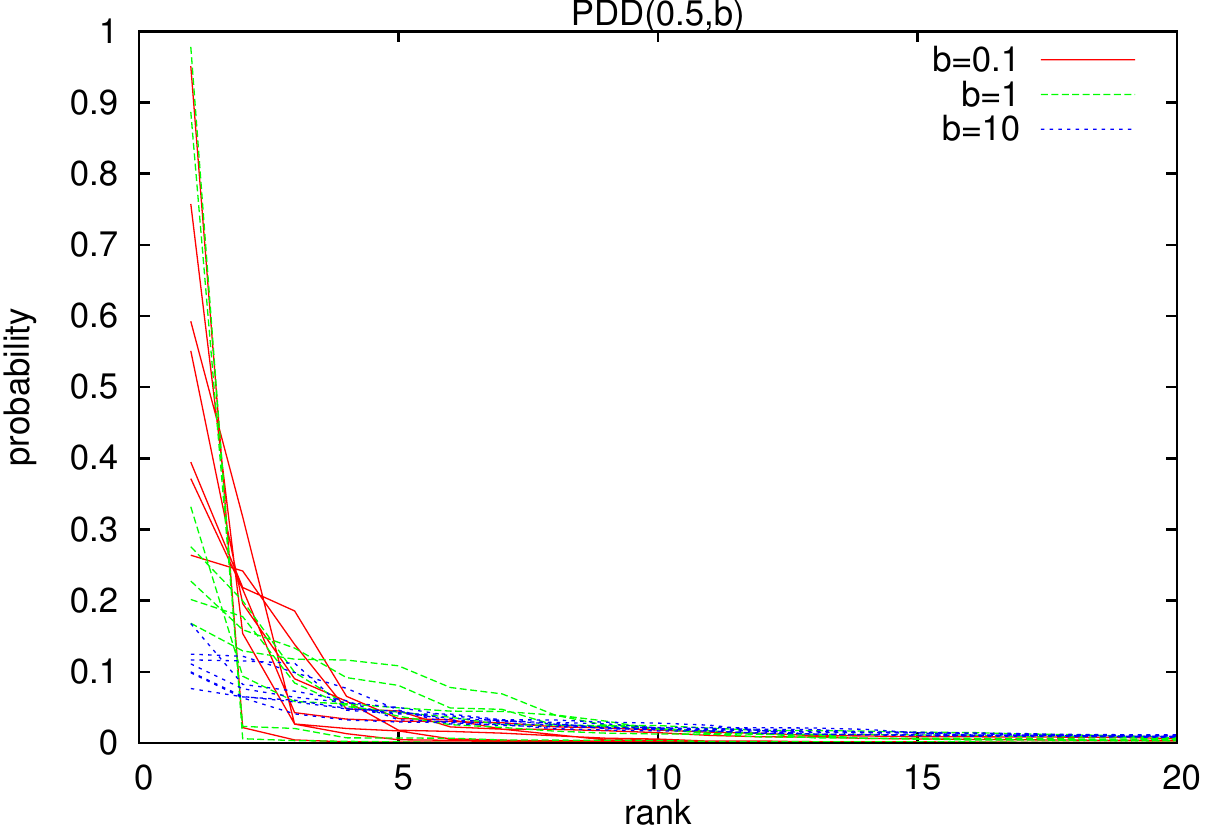}
\caption{The PDD with different parameter settings.  Each line of each
plot shows $\vec{p}\sim\mbox{PDD}(a,b)$ for a fixed discount $a$ and
a choice of concentrations $b$.  The left plots are log-log scale, showing hundreds of values.  The right plots are normal scale showing just 25 values (subsequent ones are effectively zero).}
\label{fig:pdd}
\end{center}
\end{figure*}

Why do we need two definitions, a PDD and a GEM with the
same parameters modelling different orderings of the same distribution?
\begin{itemize}
\item
When estimating the probability $\vec{p}$ for the mixture model
of Formula~\req{dfn-form} and using the GEM
distribution, the sorting is useful for sampling efficiency
\cite{KurWelTeh2007}.
One acts like one is using a PDD.
\item
The PDD gives a canonical form for the distribution.
Different sorts of the one underlying distribution can be generated
with a GEM.
\item
The GEM has a convenient sampling form
(the stick-breaking model) that allows for
simpler analysis.
\item
When used inside a mixture model or partition model,
they are indistinguishable,
so we use which ever is convenient.
\end{itemize}
Samples of $\vec{p}\sim \mbox{PDD}(a,b)$ for different parameter settings
are given in Figure~\ref{fig:pdd}, showing both probabilities and
log scale probabilities.
One can see that for the initial values (the first 20 say), the effects of
the discount $a$ are not that great, whereas the concentration
$b$ effectively changes the onset of the decrease.
We show later that vector values for the PDD with discount $a>0$
are Zipfian, behaving roughly like $p_k \propto k^{-1/a}$,
which makes the slope on the log-log plots $-1/a$.
Whereas for the case with discount $a=0$
the values are geometric, behaving roughly like $p_k \propto e^{-k/b}$.

A common definition of a Poisson-Dirichlet process is that it extends the
Poisson-Dirichlet (or GEM) distribution. This definition presents the PDP as
a functional on distributions:  it takes as input a measurable space
with domain ${\cal X}$, and a distribution over it called the {\em
base distribution},  usually represented here as $H(\cdot)$, and yields
as output a discrete distribution with a finite or countable set of
possible values on the domain ${\cal X}$.

\begin{definition}[The Poisson-Dirichlet Process]\label{defn-PDP}
Let $H(\cdot)$ be a distribution over some measurable space $\cal
X$. For $0\leq a < 1$  and $b >-a$, suppose that $\vec{p}$ is drawn
from a  Poisson-Dirichlet (or GEM) distribution with parameters $a,b$.
Moreover, let $X^*_k$ for $k=1,2,...$ be a sequence of independent
samples drawn according to $H(\cdot)$ and independent
of $\vec{p}$. Then $\vec{p}$ and $X^*_k$
for $k=1,2,...$ define a discrete distribution on $\cal X$ given by
the formula
\beq\label{dfn-form2}
  \sum_{k=1}^\infty p_k \delta_{X^*_k}(\cdot) ~.
\eeq
This distribution is a {\em Poisson-Dirichlet Process} with
parameters $a,b$ and {\em base distribution} $H(\cdot)$, denoted
$\mbox{PDP}(a,b,H(\cdot))$.
\end{definition}
The Dirichlet Process (DP) is the special case where $a=0$, and has
some quite distinct properties as shown later.
Note strict requirements on the
sorting or ordering as in Definition~\ref{defn-PY}
are not needed in the definition of the
PDP due to the effect of the summation.

The PDP is also called a {\it stochastic process} because it can be
defined as a sequence of values $X_1,X_2,... \in {\cal X}$
from some {\em base probability distribution} $H(\cdot)$ indexed by
integer valued time as $1,2,3,...$. The stochastic process is the
sequential sample from this PDP distribution. The conditional
distribution with $\vec{p}$ marginalised out for this, as long as
$H(\cdot)$ is non-atomic, is as follows:
\beq\label{eq-smp}
p(X_{N+1}\,|\,X_1,...,X_N,a,b,H(\cdot)) ~=~
  \frac{ b+Ma  }{b+N} \, H(\cdot) +
        \sum_{m=1}^M \frac{ n_m-a  }{b+N}
              \delta_{X^*_m}(\cdot) ~,
\eeq
where there are $M$ distinct values in the sequence $X_1,...,X_N$
denoted by $X^*_1,...,X^*_M$ and their occurrence counts
respectively are $n_1,...,n_M$, so $\sum_{m=1}^M n_m = N$.

The Chinese restaurant analogy for this
sequential sampling process goes as follows:
\begin{itemize}
\item A {\em customer} walks into the restaurant and sees $M$ occupied
    {\em tables} where $n_m$ others sit at table $m$ enjoying the {\em menu
    item} $X^*_m$.
\item He can start his own table with probability $\frac{ b+Ma
    }{b+N}$ and receive a new item $X^*_{M+1}$ from menu
    $H(\cdot)$ by sampling.
\item Otherwise, he goes to one of the existing $M$ tables with
  probability $\frac{ n_m-a  }{b+N}$ and enjoys the item
    $X^*_m$.
\end{itemize}
The Chinese Restaurant Process (CRP) is defined over
the partition so generated, represented with a size-biased order of
indices.
\begin{definition}[Chinese Restaurant Process]\label{defn-CRP}
For $0\leq a < 1$  and $b >-a$, generate a sequence of
integers $k_1,k_2,...$ as follows:
\[
p(k_{N+1}|k_1,...,k_N,a,b) ~=~
  \frac{ b+Ma  }{b+N} \, \delta_{k_{N+1}=M+1} +
        \sum_{m=1}^M \frac{ n_m-a  }{b+N}
              \delta_{k_{N+1}=m}~,
\]
where $M$ and the counts $n_1,n_2,...,n_M$ are derived
as for Equation~\req{eq-smp}.
\end{definition}
If indices are sampled according to the CRP
then one gets a size-biased ordered index sequence from
a $\mbox{PDD}(a,b)$ \cite{Pitman95,Ish_jasa01}, and thus the CRP
can serve as an alternative sampler for the PDP
where the probability vector $\vec{p}$ is unknown.

\section{Basic Properties}\label{sct-dd}

Before getting onto discrete domains, we review basic properties of
the PDD, and of the PDP in non-atomic domains.  Some of these
results will be used subsequently to address discrete domains.

For the sample from the distribution of Formula~\req{dfn-form2},
a latent sequence of indices
exists $I_N:=(k_1,...,k_N)$, however, these remain hidden
(only the corresponding data values are known, not the indexes).
For a non-atomic base distribution the indices are
irrelevant and we can renumber them by size-biased ordering.
Each index corresponds to a table in the CRP, and the
number of distinct indices in the sample is the number of tables
active at the restaurant.

Our notation for statistics is as follows:
\begin{definition}[Index Statistics]\label{defn-pis}
When sampling independently and identically from a discrete
distribution in the form of Formula~\req{dfn-form2}, one gets
a sequence of latent indices, an {\em index sequence} of length
$N$ given by $I_N=k_1,k_2,...,k_N$.
In $I_N$ one index value $k$ can occur multiple times.
Sort and count the $N$ points of $I_N$. Suppose there are $M$
distinct values in $I_N$, $k^*_m$ for $m=1,...,M$ that occur $n_m$
times respectively, so $\sum_{m=1}^M n_m=N$. Call $M$ the {\em
partition size} and note it depends on the sample $I_N$ and the sample size
$N$. Moreover,  use a size-biased ordering of the
indices and renumber them according to this. The
corresponding size-biased ordered index sequence is denoted $I_N^*$.
\end{definition}
Consider a sample with $N=7$ points with
latent indices $I_N=12, 435, 7198, 12, 12,  35, 7198$.
Then the size-biased ordered renumbering is $I^*_N=1,2,3,1,1,4,3$.
Thus the partition size $M=4$ and occurrence counts
$n_1=3$, $n_2=1$, $n_3=2$ and $n_4=1$.
Note the partition size $M$ for a sequence corresponds to the number
of active tables in the CRP terminology.

\begin{definition}[Data Statistics]\label{defn-pds}
When sampling from the discrete distribution of
Formula~\req{eq-smp} or via a PDP,
one gets a {\em data sequence} of length $N$
given by $S_N=X_1,X_2,...,X_N$. Sort and count the $N$ points of
$S_N$. Suppose there are $M$ distinct values in $S_N$, $X^*_m$ for
$m=1,...,M$ that occur $n_m$ times respectively, so $\sum_{m=1}^M
n_m=N$. For non-atomic base distributions $H(\cdot)$, it is safe to
associate index $m$ with data value $X^*_m$,  which is the unique
size-biased ordered renumbering. The corresponding index sequence is
denoted $I^*_N$.
\end{definition}
For discrete base distributions $H(\cdot)$, a
unique size-biased ordered renumbering for the indices does not exist
because if two data items in a sample are equal,
one cannot be certain their
latent indices are the same.

\subsection{The Dirichlet Process}\label{ssct-dp}

The Dirichlet Process (DP) is a special case of the PDP when the
discount parameter $a=0$. It has quite distinct properties
as subsequent analysis will show.  It is usually defined in
a completely different manner to the PDP as follows.
Let ${\cal X}$ be a measurable space.  For a random probability
distribution $G(\cdot)$ to be distributed according to a DP, its marginal distributions
have to be Dirichlet distributions too.
Ferguson \cite{ferguson73} gave a formal definition of the DP as follows.
\begin{definition}[Dirichlet Process]
 \label{def:DP}
   Let $H(\cdot)$ be a random measure
    on ${\cal X}$ and $b>0$ be positive
    real number. We say a random probability measure $G(\cdot)$ on ${\cal X}$
      is a Dirichlet process with a {\it base measure} $H(\cdot)$ and {\it concentration parameter}
    $b$, {\it i.e.}
    $G(\cdot) \sim \mbox{DP}(b, H(\cdot))$, if
    for any finite measurable partition $(B_1, B_2, \dots, B_k)$ of ${\cal X}$,
    the random vector $(G(B_1), G(B_2), \dots, G(B_k))$ is Dirichlet distributed
    with parameter $(b H(B_1), b H(B_2), \dots, b H(B_k))$:
  \[(G(B_1), G(B_2), \dots, G(B_k)) \sim \mbox{Dirichlet}(b H(B_1), b H(B_2),
      \dots, b H(B_k))~.\]
\end{definition}
The DP is an extension of a Dirichlet distribution, which
is defined for a finite set.
The following simple corollary of the definition above
demonstrates this.
\begin{lemma}
\label{cor:DP-DD}
According to Definition~\ref{def:DP}, if $H(\cdot)$ is a
categorical distribution over a finite space,
represented by the probability vector $\vec\theta$ say,
then the following holds
\[\mbox{DP}(b,\, H(\cdot)) = \mbox{Dirichlet}(b \vec\theta)~.\]
\end{lemma}
This is important because posterior analysis of
hierarchical Dirichlets is intrinsically difficult,
so posterior analysis of hierarchical DPs can help.

\subsection{Consistency results}\label{ssct-cons}

The PDD can be used to approximate a broader class of distributions, not
just those sampled from a $\mbox{PDD}(a,b)$. The following
lemma derived from  James \cite[Proposition~2.2]{James08} shows
this. This supposes a ``true'' probability vector $\vec{q}$ gives a
distribution of integers and then shows a sufficient property
required of $\vec{q}$ so that a PDD distribution can approximate it
based on samples.

\begin{lemma}
Suppose an integer sequence $I$ of length $N$ is sampled
independently and identically according to the probabilities
$\vec{q}$ where $0\leq q_k\leq 1$ for $k=1,2,...$ and
$\sum_{k=1}^{\infty} q_k =1$ and use the notation of
Definition~\ref{defn-pis}. If it is assumed the $\vec{q}$ is
$\mbox{PDD}(a,b)$ for $0\leq a < 1$ and $b>-a$, then the posterior
distribution on $\vec{q}$ given $I$ converges weakly to $\vec{q}$ if
$\expec{I|\vec{q},N}{M/N} \rightarrow 0$ as $N\rightarrow \infty$
where $M$ is the partition size defined in
Definition~\ref{defn-pis}.
\end{lemma}
%
Basically, we have some ``true" model over samples given by the
probability vector  $\vec{q}$. From this we compute the expected
partition size $\expec{I|\vec{q},N}{M}$ for sample sequences $I$ of
size $N$, and then check this grows slower than $N$ as $N\rightarrow
\infty$.  If this holds for $\vec{q}$, then the distribution
$\vec{q}$ can be learnt using Bayesian methods that assume $\vec{q}$
is $\mbox{PDD}(a,b)$. We show later that if $\vec{q}\sim
\mbox{PDD}(0,b)$, then almost surely $\expec{I|\vec{q},N}{M}$ is
$O(\log N)$ and if $\vec{q}\sim \mbox{PDD}(a,b)$ for $a>0$, then
almost surely $\expec{I|\vec{q},N}{M}$ is $O(N^a)$.

As a warning more than anything else, it is important to realise
the PDP should not be used to approximate continuous distributions.
This is
made precise by the consistency result due to
James \cite[Proposition~2.1]{James08}.

\begin{lemma}[PDP posterior convergence]
Suppose data is sampled independently and identically from a Polish
space $\cal X$ according to a continuous distribution $P_0(\cdot)$,
and let $H(\cdot)$ be another distribution on $\cal X$ where
$H(\cdot)$ is non-atomic.
Then the posterior of the $\mbox{PDP}(a,b,H(\cdot))$
distribution given the sample
converges weakly to point mass at the distribution
\beqn
  a H(\cdot) + (1-a) P_0(\cdot)
\eeqn
Hence the posterior is consistent only if either $H(\cdot)=
P_0(\cdot)$ or $a=0$.
\end{lemma}
Note discrete distributions cannot be continuous since they have
finite mass concentrated at points. Thus the above lemma does not
apply to the discrete case. When the ``true'' distribution
$P_0(\cdot)$ is discrete, weak converge does hold.

\subsection{Posteriors}
One can derive the probability of evidence or data given the model,
a useful diagnostic in Bayesian analysis. Various versions of this
are well known, see \cite[Appendix]{PitmanYor97} and
\cite[Proposition 9]{Pitman95}, and easily proven by induction using
the CRP.

\begin{lemma}[Probability of evidence]\label{lem-ev}
Consider finite samples $S_N=X_1,X_2,...,X_N$ from
$\mbox{PDP}(a,b,H(\cdot))$, where the base distribution $H(\cdot)$
is non-atomic. Use the notation of Definition~\ref{defn-pds}. Then
the probability of evidence given the model
$\mbox{PDP}(a,b,H(\cdot))$ is
\beqn
  p(X_1,X_2,...,X_N|a,b) =
        \frac{(b|a)_M}{(b)_N}
        \prod_{m=1}^M H(X^*_m)
   \prod_{m=1}^M (1-a)_{n_m-1}~,
\eeqn
where $(x)_N$ denotes the Pochhammer symbol
$x(x+1)...(x+N-1)=\Gamma(x+N)/\Gamma(x)$
and $(x|y)_N$ denotes
$x(x+y)...(x+(N-1)\!\cdot\! y)$, the Pochhammer symbol with increment $y$,
and $(x|0)_N=x^N$.
\end{lemma}
%

\section{Probabilities of Partitions}\label{sct-pp}

One can easily see the posterior in Lemma~\ref{lem-ev}, with the
term for the base distribution ($\prod_{m=1}^M H(X^*_m)$) removed
represents a distribution on a partition.
This section presents various properties of this
distribution.

\subsection{Chinese Restaurant Distribution}
For a partition represented as a size-biased order
$I^*_N$, the probability is a function of the
partition size $M$ and the occurrence counts
$n_1,...,n_M$, where by default the counts are listed in size-biased order.
The resultant probability on $I^*_N$ is then a neat function
$f(n_1,...,n_M)$.  This is called an
{\it exchangeable partition probability function}.
This goes under various names in the literature,
so the term {\em Chinese Restaurant Distribution} is
used here to differentiate it from the
CRP from which it is derived.
\begin{definition}[Chinese Restaurant Distribution]\label{defn-CRD}
Given a set $P$ of size $N=|P|$,
represent the partitions of $P$ by the
set of size-biased orderings of length $N$,
where one is denoted $I^*_N$.
Define
\beqn
  p(I^*_N|a,b) = \left\{
 \begin{array}{lr}
        \frac{(b|a)_M}{(b)_N}
   \prod_{m=1}^M (1-a)_{n_m-1} & \mbox{for $a>0$}  ~,\\\noalign{\medskip}
        \frac{b^M}{(b)_N}
   \prod_{m=1}^M (n_m-1)! & \mbox{for $a=0$} ~.
\end{array}
\right.
\eeqn
where $M$ and the occurrence counts $n_1,...,n_M$
follow Definition~\ref{defn-pis}.
Call this the {\em Chinese Restaurant Distribution} with
parameters $(a,b)$, abbreviated $\mbox{CRD}(P,a,b)$,
and note its samples are a partition of $P$.
\end{definition}
Note this distribution runs over all possible partitions,
so for instance if $N=3$ the possible
partitions of $\{a,b,c\}$ represented using an ordering of least
elements are given in Table~\ref{tbl-part3}.
\begin{table}
\begin{center}
\begin{tabular}{|l|ccccc|}
\hline
partition &\small$\{\{a,b,c\}\}$&\small$\{\{a,b\},\{c\}\}$&\small%
$\{\{a,c\},\{b\}\} $&\small$\{\{a\},\{b,c\}\} $&\small$\{\{a\},\{b\},\{c\}\}$\\
\hline
partition $I^*_N$ & $(1,1,1)$ & $(1,1,2)$ & $(1,2,1)$ & $(1,2,2)$ & $(1,2,3)$ \\
size $M$& 1 & 2 & 2 & 2 & 3 \\
counts $\vec{n}$ & $(3)$ & $(2,1)$ & $(2,1)$ & $(1,2)$ & $(1,1,1)$ \\
\hline
$p(I^*_N \,|\, a>0,b)$ & $\frac{(1-a)(2-a)}{(b+1)(b+2)}$ &
  $\frac{(b+a)(1-a)}{(b+1)(b+2)}$ & $\frac{(b+a)(1-a)}{(b+1)(b+2)}$ &
   $\frac{(b+a)(1-a)}{(b+1)(b+2)}$ & $\frac{(b+a)(b+2a)}{(b+1)(b+2)}$    \\
$p(I^*_N \,|\, a=0,b)$ & $\frac{2}{(b+1)(b+2)}$ &
  $\frac{b}{(b+1)(b+2)}$ & $\frac{b}{(b+1)(b+2)}$ &
   $\frac{b}{(b+1)(b+2)}$ & $\frac{b^2}{(b+1)(b+2)}$ \\
\hline
\end{tabular}
\caption{Space of partitions over $\{a,b,c\}$.}
\label{tbl-part3}
\end{center}
\end{table}

\subsection{Partition size}
The key characteristic of a sample from the
PDD or a CRD is the {\em partition size} $M$
from Definition~\ref{defn-pis}. This is related to the expected
posterior probability of starting a new bin (for the PDD or CRP) or a new
data value from $\cal X$ in the non-atomic case of the PDP.
This is given by the formula for
the unseen part of the CRP,
\beqn
  p(k_{N+1} \not\in I^*_N \,|\,  I^*_N,M, a,b)
  ~=~ p(X_{N+1} \not\in S_N \,|\, S_N,M, a,b)
  ~=~ \frac{b+M\,a}{N+b}~.
\eeqn
The posterior distribution for the partition size given just the
sample size introduces a significant function, $S^N_{M,a}$, which is
a generalised Stirling number.  It was applied to the task by Pitman
\cite[Equation~(89)]{Pitman99} where it was represented as
$a(N,M,a)$ and by Teh \cite{Teh2006a} in the form $s_a(N,M)$, as a
generalised Stirling number of type $(-1,-a,0)$ attributed to Hsu
and Shiue, where it was applied to the analysis of hierarchical
PDPs. The case for $a=0$ was first presented by Antoniak
\cite[p1161]{Antoniak74}.

\begin{lemma}[Probability on partition size]\label{lem-exp}
Consider the size-biased ordering of indices
$I_N^*$ for a sample of size $N$ from a PDD with
parameters $(a,b)$. The probability distribution for $M$ given just
$N$ and integrating over all possible partitions $I^*_N$ of
size $M$ is
\beq\label{eq-exp}
  p(M\,|\, N,a,b) ~=~ \frac{(b|a)_M}{(b)_N} S^{N}_{M,a},
  \qmbox{where}
\eeq
\beq\label{defS}
  S^N_{M,a} ~:=~ N!\,\sum_{\sum_1^M n_m = N,~n_m\ge 1}
  \prod_{m=1}^M \left( \frac{(1-a)_{n_m-1}}{n_m!}
        \frac{n_m}{\left(N-\sum_{i=1}^{m-1} n_i\right)} \right)~,
\eeq
for $M\le N$ and 0 else.
\end{lemma}
The following expressions are useful for computing $ S_{M,a}^N$.

\begin{theorem}[Expressions for $ S_{M,a}^N$]\label{thmSN} \hfill
\begin{tabbing}
\phantom{(ii)} \= \phantom{Explicit expressions} \= \kill
(i) \> Linear recursion: \>
$\displaystyle S_{M,a}^{N+1}\!\! = S_{M-1,a}^N + (N\!-\!Ma) S_{M,a}^N$
\\
    \> Boundary cond.: \>
$ S_{M,a}^N=0$ for $M>N$, $\quad S_{0,a}^N=\delta_{N,0}$.
\\
(ii) \> Mult.\ recursions:
$\displaystyle S_{M,a}^N = \sum_{n=m}^{N-M+m}{\textstyle{{N\choose n}\over{M\choose m}}} S_{m,a}^n S_{M-m,a}^{N-n}$
  $\displaystyle = \sum_{n=1}^{N-M+1}\!\!{\textstyle{N-1\choose n-1}} S_{1,a}^n S_{M-1,a}^{N-n}$
\\
  \hspace{22ex} $ S_{1,a}^N=\Gamma(N-a)/\Gamma(1-a)$. \hspace{10ex} Any $0<m<M$.
\\
(iii) \> Explicit expression: \>
$\displaystyle S_{M,a}^N = {1\over M!\,a^M}\sum_{m=0}^M{\textstyle{M\choose m}}(-)^m\prod_{h=0}^{N-1}(h\!-\!am)$ ~~~~~~for $a>0$\\
\> \> $\displaystyle
S_{M,0}^N  = {(-)^{M-1}\over (M-1)!}   \left. \frac{ \partial^{M-1}}{\partial a^{M-1} }
\left(\frac{\Gamma(N - a)}{\Gamma(1 - a)}\right) \right|^{a=0}$
\\
(iv) \> Asymptotic expr.: \>
$\displaystyle S_{M,a}^N \simeq {1\over\Gamma(1\!-\!a)}{1\over\Gamma(M)\,a^{M-1}}{\Gamma(N)\over N^a}$ for $a>0$
\\
(v) \> Expr.\ for $a=0$: \>
$\displaystyle S_{M,0}^N = |s_N^{(M)}| =$ unsigned Stirling\# of 1st kind \cite{Abramowitz:74}
\end{tabbing}
The asymptotic expression holds for $N\to\infty$ and fixed $M$ and $a$.
\end{theorem}
The explicit closed form $(iii)$ was developed in \cite{Toscano:39}
and presented in \cite{pitmanCSP}.
This explicit form becomes unstable for
large values of $M$ since it is effectively an $M$-point
interpolation to a partial derivative, thus has a lot
of differences of similar values.  It remains effective
for small $M$.  Some examples are:
\begin{eqnarray*}
S^N_{2,0} &=& S^N_{1,0} \left(\psi_0(N)-\psi_0(1)\right)
~~~~~~~~~
S^N_{3,0} = \frac{1}{2} S^N_{1,0} \left(\psi_1(N)-\psi_1(1) +
          \left(\psi_0(N)-\psi_0(1)\right)^2 \right)\\
S^N_{2,a} &=& \frac{1}{a} \left(S^N_{1,a} - S^N_{1,2a} \right)
~~~~~~~~~~~~~~~~~~~
S^N_{3,a} = \frac{1}{2a^2} \left(S^N_{1,a} - 2S^N_{1,2a}+S^N_{1,3a} \right) ~.
\end{eqnarray*}

\begin{figure*}
\begin{center}
\hspace*{-10pt}\includegraphics[width=0.52\textwidth,height=0.33\textwidth]{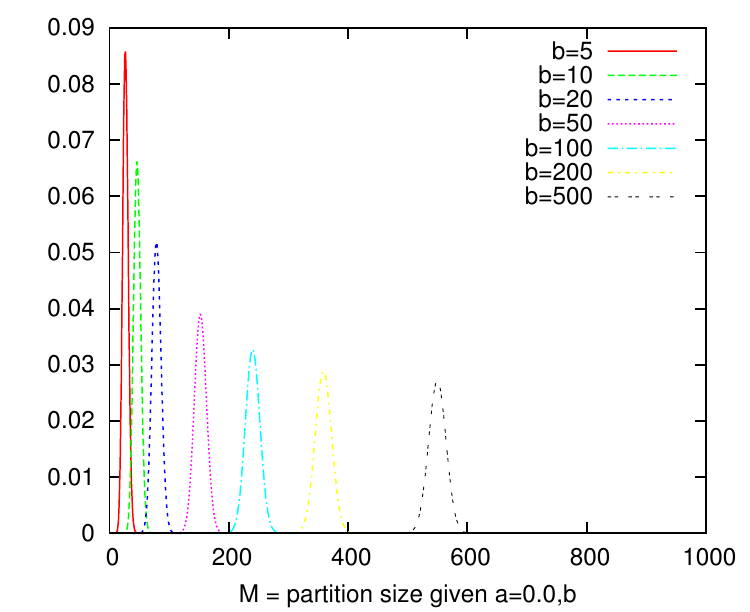}\includegraphics[width=0.52\textwidth,height=0.33\textwidth] {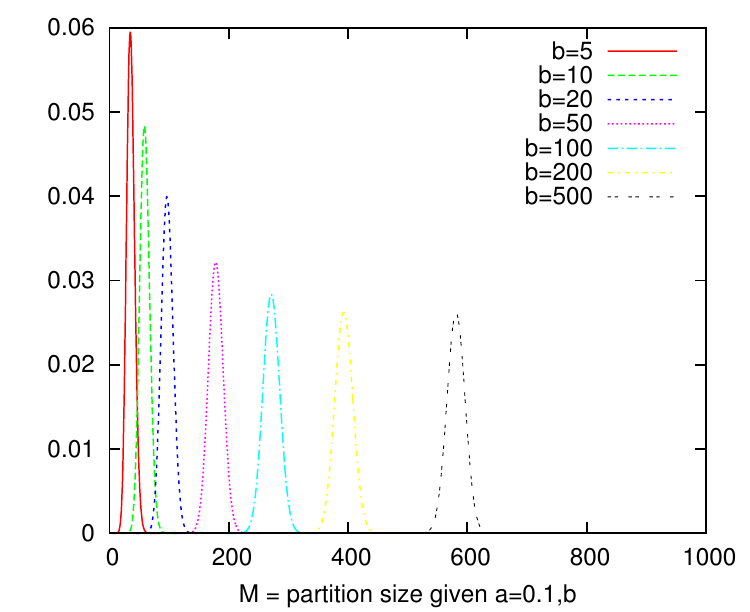}\\
\hspace*{-10pt}\includegraphics[width=0.52\textwidth,height=0.33\textwidth]{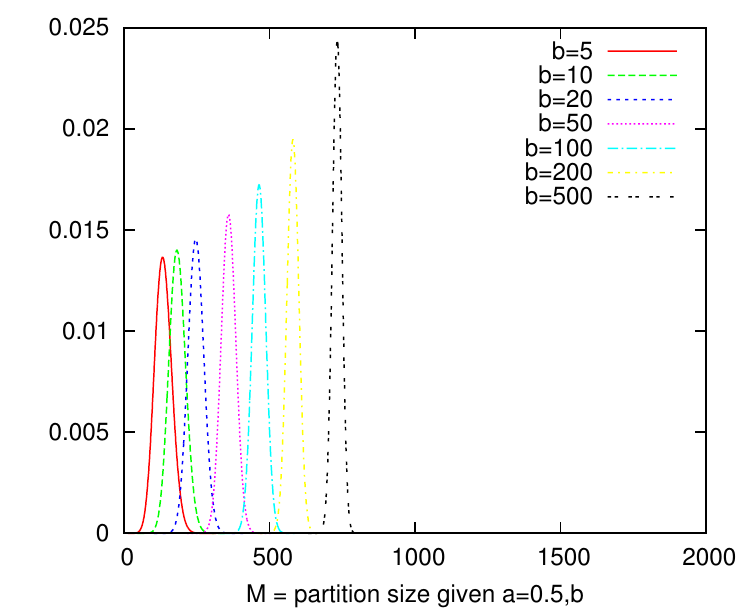}\includegraphics[width=0.52\textwidth,height=0.33\textwidth]{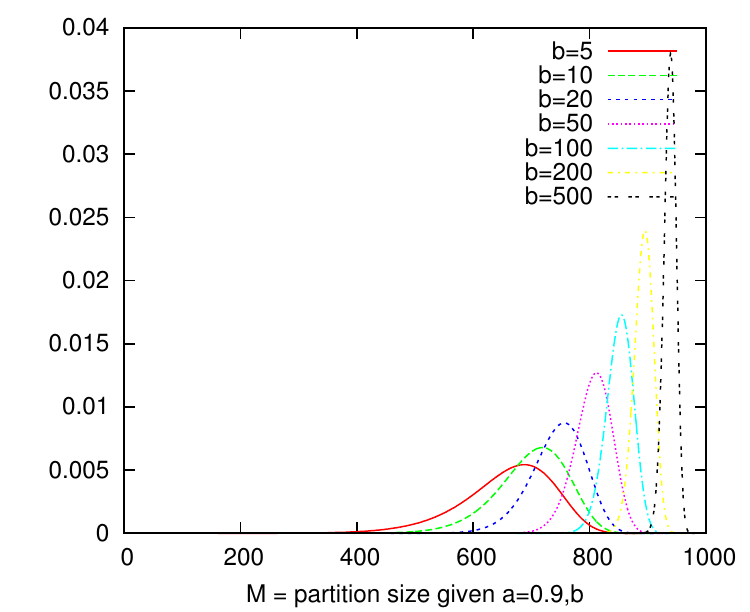}
\end{center}
\caption{Posterior probability on M given $N=1000$ and different $a$.}
\label{fig:MN}
\end{figure*}

\begin{figure*}
\hspace*{-10pt}\includegraphics[width=0.52\textwidth,height=0.33\textwidth]{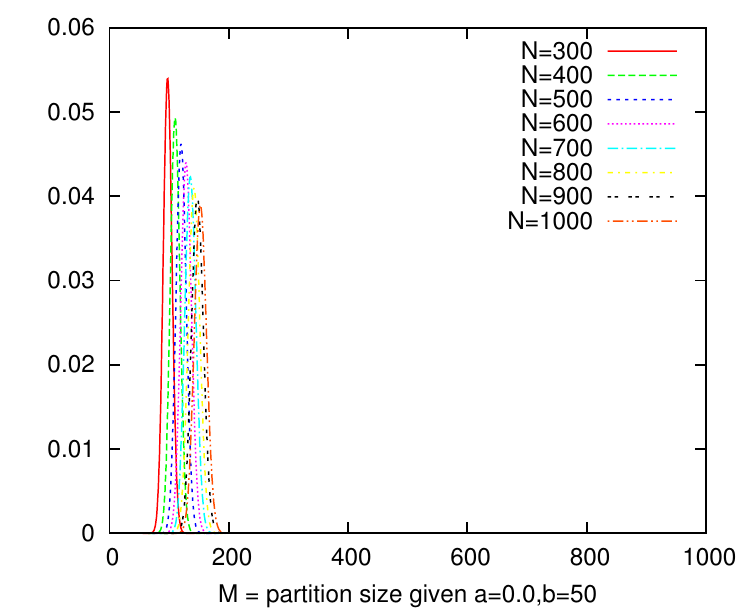}\includegraphics[width=0.52\textwidth,height=0.33\textwidth]{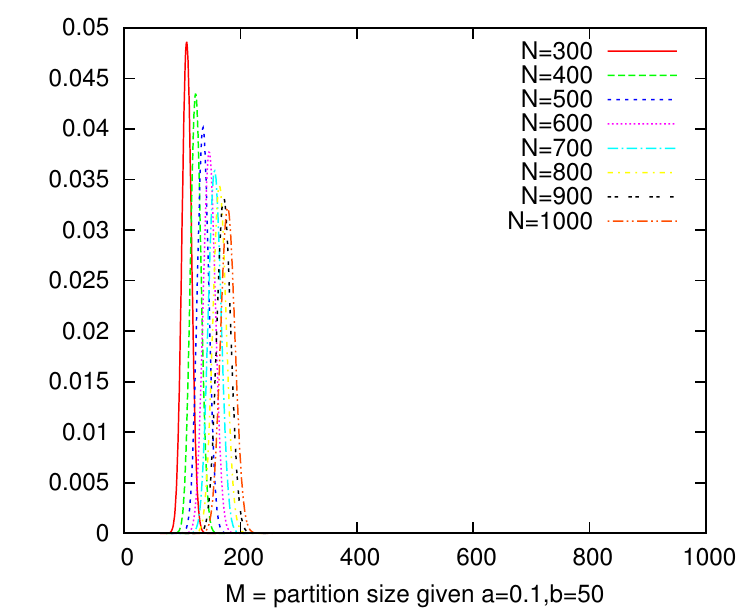}\\
\hspace*{-10pt}\includegraphics[width=0.52\textwidth,height=0.33\textwidth]{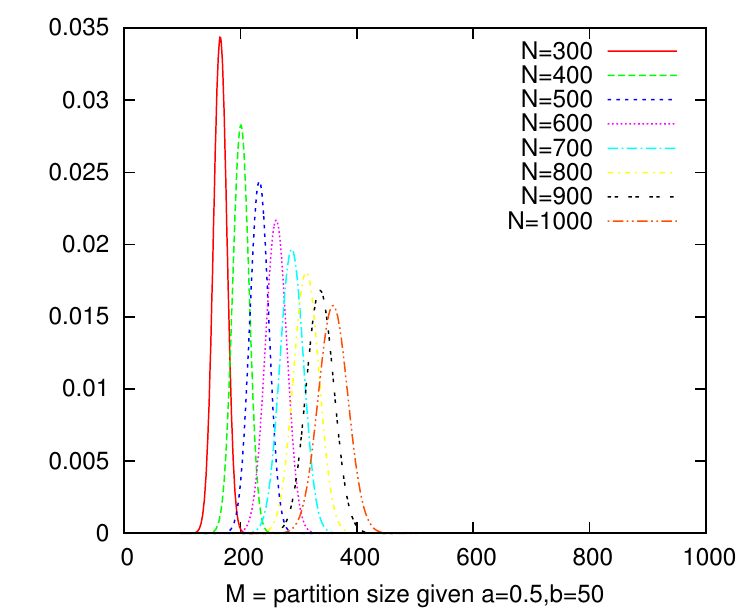}\includegraphics[width=0.52\textwidth,height=0.33\textwidth]{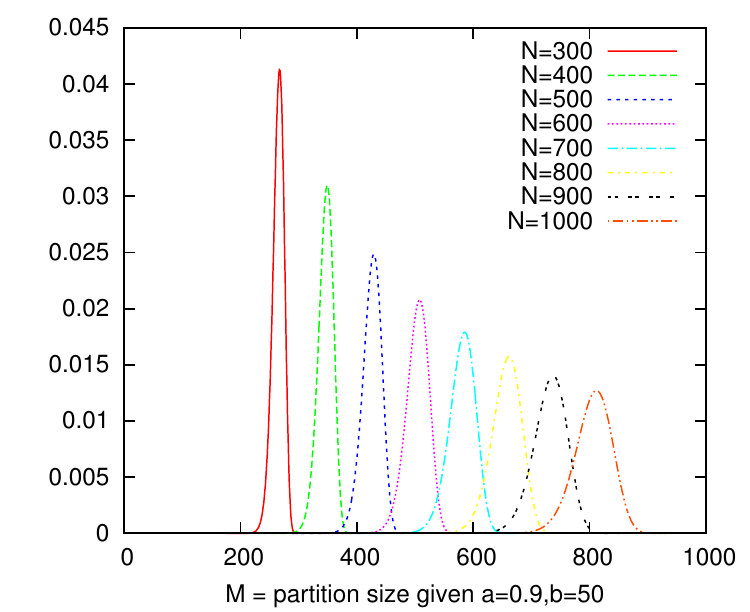}
\caption{Posterior probability on M for increasing $N$ and fixed $b=50$.}
\label{fig:bN}
\end{figure*}

Figure~\ref{fig:MN}
illustrates the shape of the distributions and their location for
different values of $a$ and $b$ and fixed $N=1000$. Similar looking
plots are produced when $N=10000$.
Note the distribution does reflect a Poisson in some ways,
being skewed both at the lower boundary $M=0$ and the upper
boundary $M=N$, and being fairly symmetric in other cases.
Figure~\ref{fig:bN} illustrates the shape of the distributions
and their location for different values of $a$ as $N$ grows,
for $b=50$.
The figure for $a=0.9$ has a different horizontal scale.
Note also how the spread of $M$ increases as the sample size $N$ increases.

\subsection{Convergence results}
It is well known that expected partition size for the DP (PDP with
$a=0$) is $O(\log N)$ and for the PDP it is
$O(N^a)$. Here the exact rates are presented
along with their expected variance \cite{YamSib00}. Further details
of moments for the PDD are also given by Ishwaran and James
\cite{Ish_jasa01}.

\begin{lemma}[Expected partition size]\label{lem-part}
In the context of Definition~\ref{defn-pis}, if a partition sequence
$I^*_N$ has
the probability vector $\vec p$ distributed {\em a priori} according
to $\mbox{PDD}(a,b)$, the expected {\em a posteriori} $M$ for a
sample of size $N$ denoted $\expec{\vec{p}|a,b,N}{M}$ (and note the
actual sample is unknown here, just its size $N$ is known), when
$a>0$ is given by
\bqan
  \expec{\vec{p}|a,b,N}{M}&=&
  \frac{b}{a}  \frac{(b+a)_N}{(b)_N}
  -  \frac{b}{a} ~,\\
  &=&
  b \left( \phi_0(b+N)-\phi_0(b)\right) + O(ab \log^2 N)\\
  &  \simeq
  &
  \frac{b}{a}  \left(1+\frac{N}{b}\right)^a  \exp\left( \frac{aN}{2b(b+N)} \right) - \frac{b}{a}
  \quad\qquad \mbox{for $N,b \gg a$} ~,\\
\eqan
where $(x)_N$ denotes the Pochhammer symbol
$x(x+1)...(x+N-1)=\Gamma(x+N)/\Gamma(x)$.
The {\em a posteriori} variance of $M$ for a sample of size $N$,
denoted $\vari{\vec{p}|a,b,N}{M}$, when $a>0$ is given by
\bqan
  \vari{\vec{p}|a,b,N}{M} &=& \frac{b(a+b)}{a^2}  \frac{(b+2a)_N}{(b)_N}
        - \frac{b}{a}  \frac{(b+a)_N}{(b)_N}
        - \left(\frac{b}{a} \frac{(b+a)_N}{(b)_N}\right)^2~ \\
  &\simeq&   \frac{b}{a}  \left(1+\frac{N}{b}\right)^{2a} \exp\left( \frac{aN}{b(b+N)} \right)
  \qquad\qquad \mbox{for $N,b \gg a$}~.
\eqan
In the context where $a=0$,
\bqan
  \expec{\vec{p}|a,b,N}{M} &=& b(\psi_0(b+N)-\psi_0(b)) \\
  &\simeq&  b\log\left(1+\frac{N}{b}\right) \qquad\qmbox{for} N,b \gg 0~,\\
  \vari{\vec{p}|a,b,N}{M}  &=& b(\psi_0(b+N)-\psi_0(b)) \\
  & & +\, b^2(\psi_1(b+N)-\psi_1(b)) \\
  &\simeq&  b\log\left(1+\frac{N}{b}\right) \qquad\qmbox{for} N > b \gg 0~,
\eqan
where $\psi_0(\cdot)$ is the digamma function and $\psi_1(\cdot)$ is
the 1-st order polygamma function, the derivative of the digamma
function.
\end{lemma}
Thus for $0\leq a < 1$ and $b$ fixed, $ \expec{\vec{p}|a,b,N}{M}$
is almost surely sublinear in $N$ as described in Section~\ref{ssct-cons}.

Note $ \expec{\vec{p}|a,b,N}{M}$ is roughly linear in $b$ in all
cases. For  the DP case (when $a=0$) and $N\gg b\gg 0$, the {\em a
posteriori} standard deviation of $M$ is approximately the square
root of $ \expec{\vec{p}|a,b,N}{M}$, so $M$ is somewhat Poisson in
its behaviour. For the PDD $a>0$ and $N\gg b\gg a$, the {\em a
posteriori} standard deviation of $M$ is approximately $
\expec{\vec{p}|a,b,N}{M}/\sqrt{b/a}$, so is smaller than $
\expec{\vec{p}|a,b,N}{M}$ for $b\gg a$.

To compare convergence of PDD distributions with known series, we
use the following lemma.

\begin{lemma}[Upper bound on expected partition size]\label{lem-series}
Suppose a partition sequence $I^*_N$ of length $N$ is sampled
independently and identically according to the probabilities
$\vec{q}$ where $0\leq q_k\leq 1$ for $k=1,2,...$ and
$\sum_{k=1}^{\infty} q_k =1$ and use the notation of
Definition~\ref{defn-pis}. If $\vec q$ takes the form of a geometric
series, $q_k=r^{k-1}(1-r)$, then
\beqn
  \expec{I^*_N|\vec q}{M} ~\leq~ \frac{\log N}{\log 1/r} + \frac{1 + 2\log 1/r+\log
\log 1/r}{\log 1/r}  ~.
\eeqn
If $\vec q$ takes the form of a Dirichlet series
$q_k=k^{-s}\zeta(s)$ for $s>1$ (where $\zeta(s)$ is the Riemann zeta
function), then
\beqn
  \expec{I^*_N|\vec q}{M} ~\leq~ 3/2 + \frac{s}{(s-1)} \left(\frac{N}{\zeta(s)}\right)^{1/s}  ~.
\eeqn
\end{lemma}
The bounds are often quite good. Experimental evaluation shows the
geometric series bound is close to about 20\% except where $r$
approaches 1, and the Dirichlet series bound is close to about 20\%
except where $s$ approaches 1.

Comparing the expected partition sizes of Lemma~\ref{lem-part} with
the different convergent series above, one can see that the PDD case
for $a>0$ behaves more like a Dirichlet series with exponent
$s=1/a$, whereas the DP case (for $a=0$) behaves more like a
geometric series with factor $r=\exp(-1/b)$.

\subsection{Dirichlet-Multinomial models}\label{sct-md}
It is instuctive to compare the
CRD with the Dirichlet-multinomial model
obtained by marginalising out parameters from a
multinomial posterior.
This is a posterior on assignments of $N$ data points
to $K$ classes, rather than a partition of $N$ data points.
So given $N$ data in total assigned to $K$ classes
with counts $\vec{n}=(n_1,...,n_K)$, and using a
prior of $\mbox{Dirichlet}_K\left(\vec{\alpha}\right)$:
\begin{eqnarray*}
p\left(\vec{n}|N,\alpha\right) &=&
\frac{\Gamma\left(b\right)}
     {\Gamma\left(N+b\right)}
\frac{\prod_k\Gamma\left(n_k+\alpha_k\right)}
     {\prod_k\Gamma\left(\alpha_k\right)}\\
&=&
\frac{1}{(b)_N}
\prod_k (\alpha_k)_{n_k}
\end{eqnarray*}
where $b=\sum_k\alpha_k$.
Here the second line represents the formula
using Pochhammer symbols to bring out the
correspondence with Definition~\ref{defn-CRD}.
Comparisons are as follows:
\begin{itemize}
\item
For the Dirichlet-multinomial, the number of classes $K$
is fixed, for the CRD the size of the partition
($M$ is sometimes used) is varied as well,
and has the posterior given in Lemma~\ref{lem-exp}.
\item
The $n_k-1$ subscript in the Pochhammer
symbol in Definition~\ref{defn-CRD} loses
1 due to starting off the new bin in the partition.
\end{itemize}

\section{Fragmentation and Coagulation}\label{sct-fc}

\begin{figure}
\begin{center}
\input{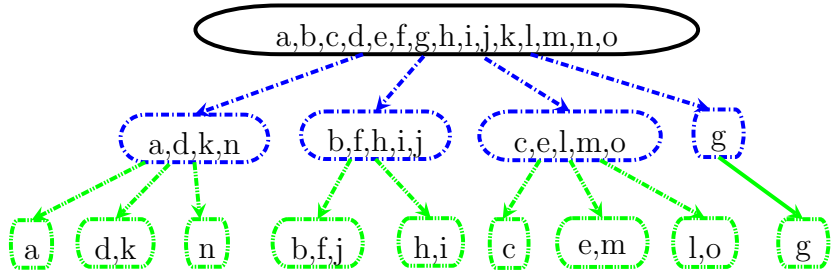}
\caption{Shows two {\em fragmentations}
 of a set of letters $\{a,b,...,o\}$.  The top node is the full set,
the middle row is after the first fragmentation, and the
bottom row is after the second fragmentation.
Each row is in size-biased order.}
\label{fig:part1}
\end{center}
\end{figure}

\begin{figure}
\begin{center}
\input{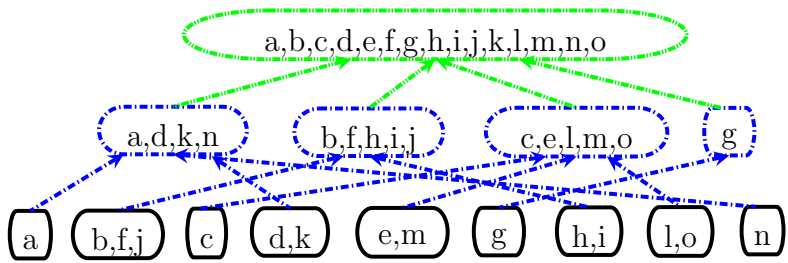}
\caption{Shows the coagulations reversing the
 fragmentations in Figure~\ref{fig:part1}.
The bottom row is the initial partition.
The middle row is a coagulation of the bottom row,
and the top row is a coagulation of the middle row.}
\label{fig:part2}
\end{center}
\end{figure}

Fragmentation is the term used when a partition
created by one distribution is further split
using partitions created by a second distribution.
In a simple finite case, the technique is illustrated in Figure~\ref{fig:part1}.
Fragmentation works as follows:
every set in the partition is further partitioned.
A complementary process is bottom-up, called coagulation, and is
illustrated in Figure~\ref{fig:part2}.
Coagulation makes a partition of the
set of sets forming the original partition,
and then merges entries in the one bin.
The use of the second partition here is a bit more indirect.
So, for instance, at the bottom row of Figure~\ref{fig:part2},
the initial set is $\{\{a\},\{b,f,j\},\{c\},\{d,k\},\{e,m\},\{g\},
  \{h,i\}, \{l,o\}, \{n\}\}$.
Partition this set into 4 parts
(1) $\{\{a\},\{d,k\}, \{n\}\}$, (2) $\{\{b,f,j\},\{h,i\}\}$,
(3) $\{\{c\},\{e,m\}, \{l,o\}\}$ and (4) $\{\{g\}\}$, and then flatten these
sets to get a partition of
$\{a,d,k,n\}\}$, $\{b,f,h,i,j\}$, $\{c,e,l,m,o\}$, and $\{g\}$
as shown in the second row.

Definitions below are given in terms of sets.  Conversion to size-biased
ordering is detailed but not difficult.
\begin{definition}[Fragmentation of partitions]\label{defn-frag}
Consider a partition $P$ represented
as a set of sets $\{P_1,P_2,...,P_M\}$,
and a sequence of partitions $Q_1,Q_2,...,Q_M$ of
the sets $P_1,P_2,...,P_M$ respectively.
Then the {\em fragmentation} of $P$ using $Q_1,Q_2,...,Q_M$
is $\bigcup_{m=1}^M Q_m$.
\end{definition}
\begin{definition}[Coagulation of partitions]\label{defn-coag}
Consider a partition $P$ represented
as a set of sets $\{P_1,P_2,...,P_M\}$,
and a second partition $Q$ of $\{1,2,...,M\}$ where $M=|P|$.
Then the {\em coagulation} of $P$ using $Q$
is $\{\bigcup_{m\in q} P_m\,:\,q\in Q\}$.
\end{definition}
As seen in the figures, fragmentation and coagulation
would seem to be complementary in their way of changing
a partition, one splits and one merges, but both driven by
partition templates.

\subsection{Operations on the CRD}
In some cases, the two operations
are also statistically complementary when applied to
samples from the Chinese Restaurant distribution.
This is usually defined in terms of fragmentation and
coagulation functionals over distributions \cite{pitmanCSP}.
For simplicity we present the operations
directly in terms of Definitions~\ref{defn-frag} and~\ref{defn-coag}
as follows:
\begin{theorem}[Fragmentation of CRDs]\label{th-pfrag}
For $0< a_1 < 1$, $0\leq a_2 < 1$  and $b >-a_1a_2$
and a set $P$,
if a partition $Q$ sampled from $\mbox{CRD}\left(P,a_1a_2, b\right)$
is fragmented with
partitions sampled from  $\mbox{CRD}\left(Q_m, a_1, -a_1a_2\right)$,
for $Q_m\in Q$,
then the resultant fragmented partition
is distributed as $\mbox{CRD}\left(P, a_1, b\right)$.
\end{theorem}
\begin{theorem}[Coagulation of CRDs]\label{th-pcoag}
For $0< a_1 < 1$, $0\leq a_2 < 1$  and $b >-a_1a_2$ and a set $P$,
let $I \sim \mbox{CRD}\left(P, a_1, b\right)$,
the partition size of $I$ is $M$,
and $J \sim \mbox{CRD}\left(I, a_2, b/a_1\right)$
then a coagulation of $I$ with $J$
is $\mbox{CRD}\left(P, a_1a_2, b\right)$.
\end{theorem}
Thus one coagulates to convert a partition from
$\mbox{CRD}\left(P, a_1, b\right)$ to be a partition
from $\mbox{CRD}\left(P, a_1a_2, b\right)$,
and conversely, one fragments to convert a partition
from $\mbox{CRD}\left(P, a_1a_2, b\right)$
to be a partition
from $\mbox{CRD}\left(P, a_1, b\right)$.

Note that if we traverse down a tree
levelwise, it corresponds to performing a sequence of
fragmentations.
Likewise, if all leaves occur at the same depth
(a condition that can be enforced by the insertion
 of dummy nodes), then one can traverse up the tree
levelwise, and it corresponds to performing a sequence of
coagulations.
If we have a schedule of strictly increasing discounts
$a_1,a_2,a_3, ...$ and a maximum tree depth, then we can generate a tree
with leaves in $\{1,2,...,N\}$ levelwise as given in Algorithm~\ref{alg:tree}.
\begin{algorithm}
  \algsetup{indent=2em, linenodelimiter=.}
  \centering
  \begin{algorithmic}[1]
    \STATE $root=\{1,2,...,N\}$;
    \STATE $tree=()$;    $nodes=()$;
    \STATE $list \sim \mbox{CRD}(root,a_{1},b)$;
    \FOR{ $m\in list$ }
       \STATE $push(tree,parent(m,root)) $;
       \STATE $push(nodes,(m,1)) $;
    \ENDFOR
    \WHILE{$(node,depth)=pop(nodes)$}
      \IF { $depth<maxdepth$ }
        \IF { $|node|>1 $ }
              \STATE  $list \sim \mbox{CRD}(node,a_{depth+1},-a_{depth})$;
              \FOR{ $m\in list$ }
                 \STATE $push(tree,parent(m,node)) $;
                 \STATE $push(nodes,(m,depth+1)) $;
              \ENDFOR
        \ELSE
            \STATE \COMMENT{Copy through to next depth.}
            \STATE $push(tree,(node,depth+1)) $;
            \STATE $push(nodes,(node,depth+1)) $;
        \ENDIF
      \ENDIF
    \ENDWHILE
    \STATE \COMMENT{Tree is represented as a list of $parent(\cdot,\cdot)$
      relations.}
    \RETURN $tree$;
  \end{algorithmic}
\caption{Sampling a tree with $N$ nodes using schedule $a_1,a_2,...,a_{maxdepth}$.}
\label{alg:tree}
\end{algorithm}
For the partitions generated, it is readily seen that
\begin{itemize}
\item
 the partition at depth $D$
(such that $1\leq D\leq maxdepth$) as generated
will be distributed as $\mbox{CRD}\left(root,a_{D},b\right)$;
\item
the partition applied to a single node
(excepting the root node) at level $D$ is
$\mbox{CRD}\left(node,a_{D+1},-a_{D}\right)$.
\end{itemize}
The expected partition sizes can then be inferred from
Lemma~\ref{lem-part}.
For instance, if $a_1=0$, then
the first level below the root will have on average
approximately $b\log\left(1+N/b\right)$ nodes.
Subsequent node partitioning will occur independently
of $b$ as the subsequent CRD parameters are set
by the discount schedule only.  This is illustrated
\begin{figure*}
\includegraphics[width=0.49\textwidth,height=0.4\textwidth]{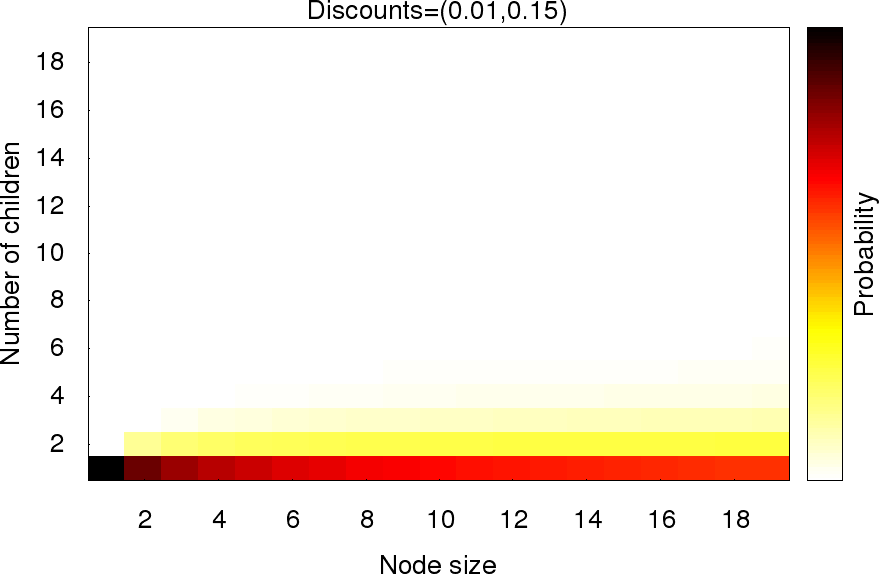}\hspace*{12pt}\includegraphics[width=0.49\textwidth,height=0.4\textwidth]{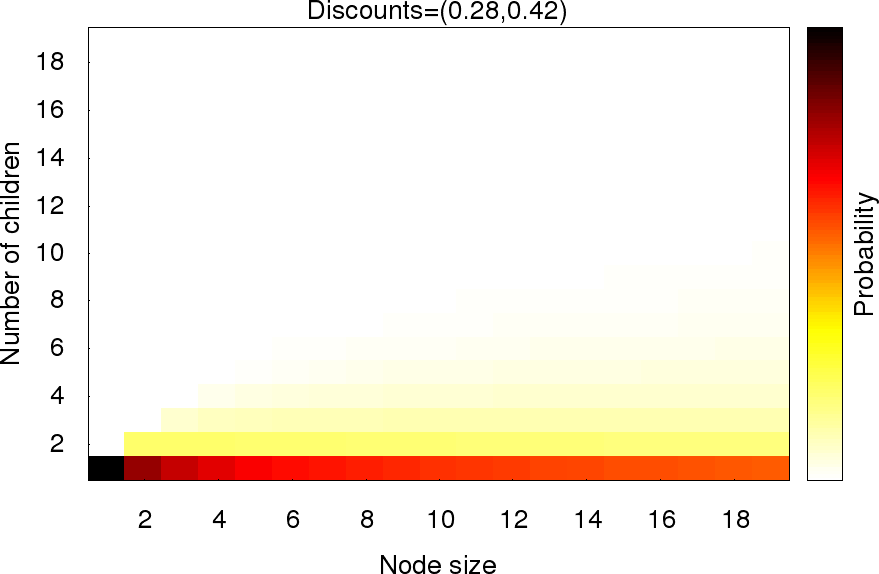}\\~\\
\includegraphics[width=0.49\textwidth,height=0.4\textwidth]{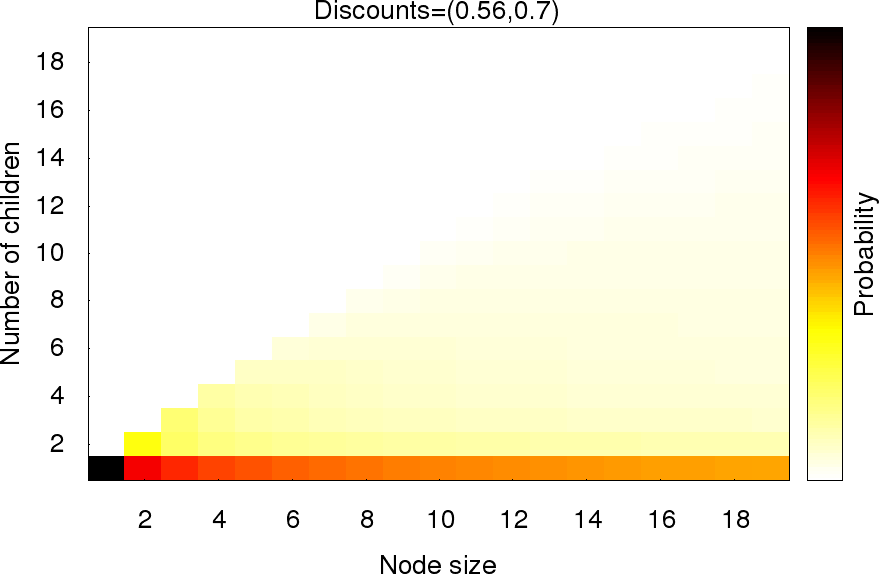}\hspace*{12pt}\includegraphics[width=0.49\textwidth,height=0.4\textwidth]{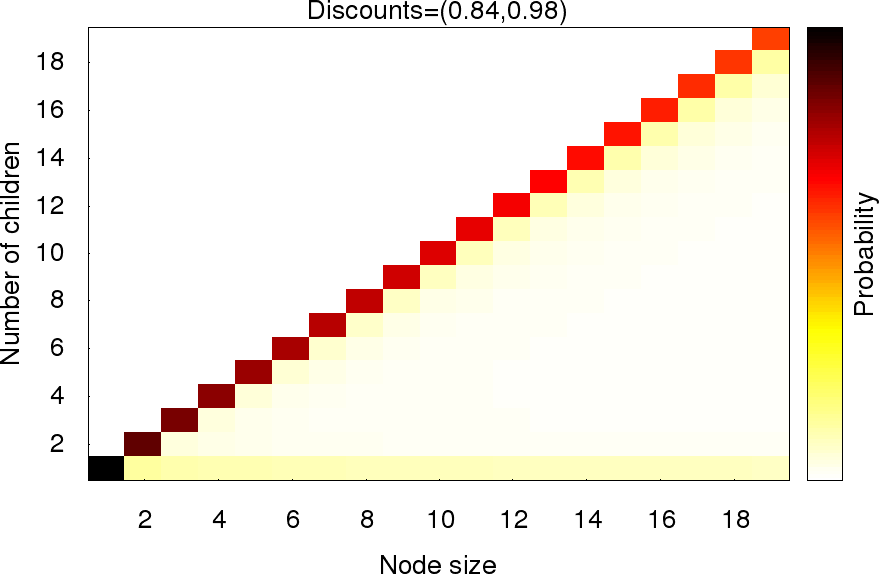}
\caption{Distribution over the number of children for internal nodes based on discount schedule, thus using $\mbox{CRD}\left(node,a_{D+1},-a_{D}\right)$.}
\label{fig:tree}
\end{figure*}
in Figure~\ref{fig:tree}.
For instance, when splitting at a level where the discount is $0.56$ and
the discount for the next level down is $0.70$, then look at the
third plot with discounts labelled $(0.56,0.70)$
which represents applying a $\mbox{CRD}(node,0.70,-0.56)$.
By the orange colour,
about 50\% of nodes of size 10 remain unpartitioned, but the
remained are split into 2-8 say children.
Whereas for the fourth plot with discounts labelled $(0.84,0.98)$,
80\% of nodes are maximally partitioned into single element nodes
(leaves) and some small percentage will remain unpartitioned.
One can see from this that by judicious use of a discount schedule, one
can generate trees that split with increasing propensities at
lower levels.

\subsection{Operations on the PDD and PDP}
Fragmentation and coagulation also apply to
the distributions on infinite partitions
represented by PDDs (or GEMs) and their
counterpart the PDP.

The basic construction is most easily understood in terms
of the basic summation form of Formula~\req{dfn-form2}.
We start with a simple draw from a PDP,
$\sum_{k=1}^\infty p_k \delta_{X^*_k}(\cdot)$.
However, we now replace the draws from the base
distribution $X^*_k$ by another draw from a PDP
with probability vector $\vec{q}_k$, so we get
\begin{equation}
\label{eq-frag}
\sum_{k=1}^\infty p_k \sum_{j=1}^\infty q_{k,j} \delta_{X^*_{k,j}}(\cdot)
~=~
\sum_{k,j=1}^\infty (p_kq_{k,j}) \delta_{X^*_{k,j}}(\cdot)
~.
\end{equation}
In some circumstances the terms $p_k q_{k,j}$ can be shown
to follow a PDD distribution.
The general result \cite{pitmanCSP} reworded for simplicity is as follows:
\begin{theorem}[Fragmentation of PDDs]\label{th-frag}
For $0< a_1 < 1$, $0\leq a_2 < 1$  and $b >-a_1a_2$,
if $\vec{p}\sim \mbox{PDD}\left( a_1a_2, b\right)$
and $\vec{q}_k \sim \mbox{PDD}\left( a_1, -a_1a_2\right)$
for each $k=1,...,\infty$
then a sort of the resultant
$\{p_kq_{k,j}\,:\,j,k=1,...,\infty\}$
is  $\mbox{PDD}\left( a_1, b\right)$.
\end{theorem}
\begin{corollary}[Fragmentation of PDPs]
Let $H(\cdot)$ be a distribution over some measurable space $\cal X$.
For $0< a_1 < 1$, $0\leq a_2 < 1$  and $b >-a_1a_2$
introduce a latent distribution $Q(\cdot)$
\[
Q(\cdot)~\sim~
\mbox{PDP}\left( a_1a_2, b, \mbox{PDP}\left( a_1, -a_1a_2, H(\cdot)\right)\right)
~.
\]
A sample from $R(\cdot)$ is taken by first drawing a
sample $Q(\cdot)$ by the nested PDP above,
and then taking a sample from $Q(\cdot)$.
Then it follows that
\[
R(\cdot) ~\sim~\mbox{PDP}\left( a_1, b, H(\cdot)\right) ~.
\]
\end{corollary}

Coagulation is the inverse of fragmentation so
one follows a similar explanation.
Start with a simple draw from a PDP,
$\sum_{k=1}^\infty p_k \delta_{Y^*_k}(\cdot)$.
However, we now wish to cluster some of the
values $\delta_{Y^*_k}(\cdot)$.  To do this we obtain a
second draw from another PDP,
$Q(\cdot)=\sum_{j=1}^\infty q_j \delta_{X^*_j}(\cdot)$,
and replace each $Y^*_k$ by a sample $X^*_{j_k}\sim Q(\cdot)$.
So we get
\begin{equation}
\label{eq-coag}
\sum_{k=1}^\infty p_k  \delta_{X^*_{j_k}}(\cdot)
~=~
\sum_{j=1}^\infty \left( \sum_{k\,:\,j_k=j}p_k \right) \delta_{X^*_{j}}(\cdot)
~.
\end{equation}
In some circumstances the terms $\sum_{k\,:\,j_k=j}p_k$ can be shown
to follow a PDD distribution.
The general result \cite{pitmanCSP} again is as follows:
\begin{theorem}[Coagulation of PDDs]\label{th-coag}
For $0< a_1 < 1$, $0\leq a_2 < 1$  and $b >-a_1a_2$,
if $\vec{p}\sim \mbox{PDD}\left( a_1, b\right)$
and $\vec{q} \sim \mbox{PDD}\left( a_2, b/a_1\right)$
and for $k=1,...,\infty$, $j_k\sim \vec{q}$,
then a sort of the resultant
$\left\{ \sum_{k\,:\,j_k=j} p_k\,:\,j=1,...,\infty\right\}$
is  $\mbox{PDD}\left( a_1a_2, b\right)$.
\end{theorem}
\begin{corollary}[Coagulation of PDPs]
Let $H(\cdot)$ be a distribution over some measurable space $\cal X$.
For $0< a_1 < 1$, $0\leq a_2 < 1$  and $b >-a_1a_2$,
sample a distribution $R(\cdot)$ using
a latent distribution $Q(\cdot)$ as follows:
\[
Q(\cdot)~\sim~\mbox{PDP}\left( a_2, b/a_1, H(\cdot)\right)
~~~~~~~~~~~~~~~~~~~
R(\cdot)~\sim~
\mbox{PDP}\left( a_1, b, Q(\cdot)\right)~.
\]
Marginalising out $Q(\cdot)$, it follows that
\[
R(\cdot) ~\sim~\mbox{PDP}\left( a_1a_2, b, H(\cdot)\right) ~.
\]
\end{corollary}
The intrinsic nature of fragmentation and coagulation is
typified by Equations~\req{eq-frag} and~\req{eq-coag}
respectively and they demonstrate a duality:
for $0\leq a_2<a_1<1$ and $b>=-a_2$,
\[
\mbox{PDD}(a_1,b)
~\begin{array}{c}
\stackrel{fragmentation}{\Longrightarrow}\\
\\
\stackrel{coagulation}{\Longleftarrow}
\end{array}
~\mbox{PDD}(a_2,b)
\]
Moreover, fragmentation is an operator for simplifying
nested PDPs whereas coagulation is for simplifying
hierarchical PDPs
(as used in \cite{WooArcGas2009a}).

\section{Improper Priors}\label{sct-dm}

A distribution or prior is called {\it proper} if it integrates (or
sums) to one. The Bayesian theory of {\it improper priors} allows
one to extend the space of reasonable priors.  The idea is that if
the posteriors from the prior are always proper, then perhaps one
can represent the improper prior as a sequence of proper priors. The
limit of this sequence may not be proper, but at least its
posteriors all are. In  this section we develop an improper prior
that corresponds to the PDD.

With any $L_d$ distance for $d\geq 1$, the infinite-dimensional
probability vector $\vec p$ of Formula~\req{dfn-form2} defines a
Hilbert space\footnote{Only when $d\geq 1$ is the subsequent
distance guaranteed to be finite for any two members of the space.}.
It is difficult to define a prior probability on such a space
because not only does one require a measure be defined for the
infinite vector, it must be normalised, so the total measure is 1.
Some theories just give priors for finite linear projections of the
full Hilbert space, for instance the cylindrical measures of Minlos
\cite{minlos01}. This is sufficient according to Carathodory's
extension theorem, see Bogachev \cite{Bogachev:07}, to define the
prior on the full space\footnote{The cylinders form a semi-ring, and
we have a countably additive (pre-)measure on the semi-ring, this
implies a unique extension on the generated ring, the sigma-algebra
is generated by the cylinder sets, and Caratheodory's extension
theorem shows that there exists a unique extension of the
(pre-)measure to the sigma-algebra.}. For the PDD model, an
additional problem is the projections of the prior on finite vector
subspaces appear to be improper as well. Thus, the best one can do
is define a prior in terms of a measure for all finite sub-vectors
as follows:

\begin{definition}[Improper prior for PDDs]\label{defn-pyd}
Given parameters $(a,b)$, where $0\leq a <1$ and $b>-a$, define the
{\em improper prior for PDDs} (an unnormalised measure) as follows.
Take any reordering of the infinite-dimensional probability vector
$\vec p$, and then for every sub-vector $p_1,p_2,...,p_M$ of the
reordering, use the  following measure:
\beqn
  p(p_1,p_2,...,p_M,p^+_M) ~:=~
  \left(p^+_M\right)^{b+M\,a-1} \prod_{m=1}^M p_m^{-a-1}  ~,
\eeqn
where $p^+_M=1-\sum_{m=1}^M p_m$.
\end{definition}
Note this applies to every sub-vector, so ordering of the
probabilities is not needed as in Definition~\ref{defn-PY}. The
measure $p(p_1,p_2,...,p_M,p^+_M)$ in the definition is an instance
of an $M+1$-dimensional improper Dirichlet with parameters
$(-a,-a,...,-a,b+M\,a)$, denoted here informally as
\beqn
  \mbox{Dirichlet}_{M+1}(-a,-a,...,-a,b+M\,a) ~.
\eeqn
Moreover, note that we believe this measure has no corresponding
limit form as $M\rightarrow \infty$ on the full infinite-dimensional
probability vector $\vec p$. Given an improper prior measure, one
can infer a posterior measure using an unnormalised version of Bayes
theorem.  If the posterior measure can be normalised, then the
posterior is now a correct probability.

It is shown next that the definition is consistent in the sense that
the measures for different sub-vectors are natural extensions of one
another. This property is called additivity for proper Dirichlets
and is well-known. It is plausible that it should hold for the
improper measure too, but the standard proofs cannot be transferred
since the involved integrals no longer exist.  Here we check
additivity does transfer to improper Dirichlets.

\begin{lemma}[Consistency of projections]\label{lem-cons}
In the context of Definition~\ref{defn-pyd}, if the prior measure
for $p_1,...,p_M$ is projected down to some sub-vector, say
$p_1,...,p_{L}$ for $L<M$, then the projected measure is consistent
with Definition~\ref{defn-pyd}.
\end{lemma}
We must now show the improper prior for PDDs is well defined. This
is done using the $L_{1}$ (total variation) distance defined for
probability density functions $H(\cdot)$ and $G(\cdot)$ as follows
\begin{equation}
  L_{1}(H,G) ~=~
   \int_{\vec{p}}
       \left| H(\vec{p}) - G(\vec{p}) \right| \mbox{d}\,\vec{p} ~.
\label{eqn-lh}
\end{equation}
The theorem below says that a sequence of proper priors exist that
can approximate the  improper prior for PDDs arbitrarily closely in
the sense that their posteriors given any sequence sample can be made
arbitrarily close to the corresponding proper posterior of the
improper prior. Closeness here is measured by total variation
distance.

\begin{theorem}[Justifying improper prior]\label{thm-dap}
Using the notation of Definitions~\ref{defn-pyd} and~\ref{defn-pis},
there exists a set of proper priors $G_{\delta}$ for $\delta>0$ such
that for  any $\epsilon>0$ and any sample $I_N$ there exists a
$\delta_0$ such that for all $0<\delta<\delta_0$ the proper
posterior~(I) given  $I^*_N$ of Lemma~\ref{lem-base} is within
$\epsilon$ by the $L_{1}$ distance of the posterior of $G_{\delta}$
given  $I^*_N$.
\end{theorem}

Because the improper prior is well defined, one can justifiably
obtain posteriors and sampling results from the prior. Now these are
identical to those for the PDD and PDP, as we detail below, however
they were derived from the improper prior, not from any of the
standard definitions for PDDs or PDPs.

\begin{lemma}[Some properties]\label{lem-base}
Using the improper prior for PDDs with parameters $(a,b)$ and
non-atomic base distribution $H(\cdot)$, the following holds:
\begin{description}
\item[Proper posteriors (I):]
Using the notation of Definitions~\ref{defn-pis} and~\ref{defn-pyd}.
The posterior distribution given $I^*_N$ is
\beq\label{eq-pddp}
   \left(p_1,...,p_M,p^+_M\right)  \,|\, I^*_N
       ~\sim~ \mbox{Dirichlet}(n_1-a,...,n_M-a,b+Ma)~,
\eeq
where $p^+_M=1-\sum_{m=1}^M p_m$.
\item[Proper posteriors (II):]
Using the notation of Definition~\ref{defn-pds},
in the case of arbitrary samples $S_N$, the
posterior distribution given $S_N$ is
\bqa\label{eq-pdpp}
  X_{N+1}\,|\, p_1,...,p_M,p^+_M, S_N
  & \sim & p^+_M \, H(\cdot) + \sum_{m=1}^M p_m \delta_{X^*_m}(\cdot) \\ \nonumber
  \left(p_1,...,p_M,p^+_M\right)  \,|\, S_N
  & \sim & \mbox{Dirichlet}(n_1-a,...,n_M-a,b+Ma)~,
\eqa
where $p^+_M=1-\sum_{m=1}^M p_m$.
\item[Sampling: ]
Using the notation of Definition~\ref{defn-pds},
if we marginalise out the probability vector $\vec{p}$, then
the posterior distribution in the next sample $X_{N+1}$,
$p(X_{N+1}\,|\,S_N)$, is
\beqn
  \frac{ b+Ma  }{b+N} \, H(\cdot) +
        \sum_{m=1}^M \frac{ n_m-a  }{b+N}
              \delta_{X^*_m}(\cdot) ~.
\eeqn
\item[Stick-breaking: ]
A stick-breaking like construction holds for the posteriors~(I) and~(II) above.
That is, for $1\leq m\leq M$
\beqn
  p_m~=~ V_m \prod_{i=1}^{m-1}(1-V_i)
\eeqn
where each $V_m$ is
independent $\mbox{Beta}\left(n_m-a,b+ma+\sum_{i=m+1}^M n_i\right)$.
Since the $\left(n_m,\sum_{i=m+1}^M n_i \right)$ are count terms,
one can say each $V_m$ has an improper prior $\mbox{Beta}(-a,b+ma)$.
\item[Size-biased sampling: ]
The prior on a size-biased ordering from
the improper prior for PDDs with parameters $(a,b)$
follows a $\mbox{GEM}(a,b)$.
\end{description}
\end{lemma}
The posterior formulation for PDPs corresponding to
Equation~\req{eq-pdpp} is attributed to Pitman \cite{Pitman96} by
Ishwaran and James  \cite[Section~4.4]{Ish_jasa01}. The sampling
result is the standard Chinese Restaurant Process for the PDP from
Ishwaran and James  \cite[Section~2.2]{Ish_jasa01}. The
stick-breaking result here is different to the standard PDP
\cite[Section~2.1]{Ish_jasa01}, which has stick priors
$\mbox{Beta}(1-a,b+ma)$ (that is, it is proper), see
\cite{PitmanYor97}. Here we use improper priors
$\mbox{Beta}(-a,b+ma)$, which matches the sampling of the CRP as
described above.

Now the $\mbox{PDD}(a,b)$ distribution is defined with sorting,
whereas the improper prior for PDDs with parameters $a,b$ is not.
Therefore they do not correspond directly unless some sorting is
done. So we can say that sorting the $\vec{p}$ for an
improper prior for PDDs yields a $\mbox{PDD}(a,b)$ distribution.

\section{The Discrete Case}\label{sct-post}

Now consider the case of discrete base distributions.
In this case, $H(\cdot)$ is a probability function,
not a probability density function, so for each sample $X_k$  from
$H(\cdot)$ its probability is finite and thus identical draws can be
repeated.
This makes the evidence calculation of Lemma~\ref{lem-ev}
invalid whenever the PDP is used in discrete or hierarchical contexts.
Here we present the techniques used to get around this.

If we have not
been given the index sequences for a sample from
an impulse mixture model,
we can only guess what the
indices might be. In this case, the
detail of the data partition is partially
hidden. So for instance, consider the sample of words:
\begin{quote}
  ``from",``apple",``to",``from",``from",``cat",``to",...
\end{quote}
This can have a size-biased ordered index sequence of
$I^*_N=1, 2, 3, 1, 1, 4, 3$.
However, since $H(\cdot)$ is discrete, it could be that
the three instances of  ``from" come from different indices
in Formula~\req{dfn-form2}.
Some other size-biased orderings of indices compatible with this
sequence of words are as follows:
\begin{quote}
  1, 2, 3, 1, 1, 4, 3, ...,~~~~~~~1, 2, 3, 1, 4, 5, 3, ...,\\
  1, 2, 3, 1, 1, 4, 5, ...,~~~~~~~1, 2, 3, 4, 5, 6, 3, ...,
\end{quote}
We are unable to say which is correct, however, they
have a single coarsest version.
Thus we introduce additional latent variables to account
for the uncertainty, introduced next.

\subsection{Multiplicity}
The definition below, {\it
multiplicity}, measures the cardinality of the (unknown) set of
indices contributing to one observation $X$ that occurs multiple times
in the data.  In the Chinese restaurant analogy, the multiplicity
for the data value $X$
is the number of active tables  with the particular menu item $X$.

\begin{definition}[Multiplicity]\label{dfn-sm}
Consider Definition~\ref{defn-pds}, but now
assume that the base distribution $H(\cdot)$ is discrete.
For a given sample $S_N$,
let $I_N$ be a size-biased ordered index sequence matching $S_N$,
and assume they are represented as
$S_N=(X_1,...,X_N)$ and $I_N=(k_1,...,k_N)$.
The {\em multiplicity} of the value $X \in S_N$ is defined as the size of
the set $\{ k_n\,:\, n=1,...,N, \, X_n=X\,\}$.
\end{definition}
Multiplicities are statistics from the latent indices $I_N$ and are
thus themselves latent.
Continuing the words example at the start of Section~\ref{sct-post},
some of the potential size-biased ordered index sequences
are illustrated in Table~\ref{tbl-mult}.
The first column gives the map from index to word,
the second column gives the resultant index sequence,
and the third column lists the multiplicities for the
words ``from",``apple",``to",``cat" respectively
({\it i.e.}, in size-biased order).
Note, for instance, in the last row, the word
``from'' appears 3 times in the map and thus has
multiplicity 3.
\begin{table}
\begin{center}
\begin{tabular}{|l|l|c|}
\hline
Word map & Index sequence $I_N$ & Multiplicities $\vec{t}$ \\\hline
``from",``apple",``to",``cat" &1, 2, 3, 1, 1, 4, 3 & 1,1,1,1\\
``from",``apple",``to",``from",``cat" &1, 2, 3, 1, 4, 5, 3 & 2,1,1,1\\
``from",``apple",``to",``cat",``to"&1, 2, 3, 1, 1, 4, 5 & 1,2,1,1\\
``from",``apple",``to",``from",``from",``cat"&1, 2, 3, 4, 5, 6, 3 & 3,1,1,1\\
\hline
\end{tabular}
\caption{Multiplicities from different index sequences of the sequence
  $S_7=\{ \mbox{``from",``apple",``to",``from",``from",``cat",``to''}\}$.}
\label{tbl-mult}
\end{center}
\end{table}

For the discrete base distribution,
we must consider the situation where the
multiplicities can be greater than one, so a more general
probability of evidence result is needed, for instance for
Lemma~\ref{lem-ev}, since $\mbox{PDP}(a,b,H)$ returns values from
$\cal X$, but no indices. The following corollary of
Lemmas~\ref{lem-ev} and~\ref{lem-exp} is a special case of
\cite[Equation~(31)]{Teh2006a}, there proven directly for the
hierarchical PDP.

\begin{corollary}[Evidence for discrete case]\label{lem-ev2}
Consider the probability of evidence for a finite sample
$X_1,X_2,...,X_N$ from $\mbox{PDP}(a,b,H)$ with discrete base
distribution $H(\cdot)$. Use Definition~\ref{defn-pds}, and let
$t_m$ be the latent multiplicity of $X^*_m$ in the sample, and let
their total $\sum_{m=1}^M t_m = T$. Note they must satisfy the
constraints $0\leq t_m\leq n_m$ and $t_m=0$ if and only if $n_m=0$.
Then the joint probability of the sample and the multiplicities is:
\beqn
  p(X_1,X_2,...,X_N,t_1,...,t_M\,|\,a,b,H(\cdot)) ~=~
    \frac{(b|a)_{T}}{(b)_N}
   \prod_{m=1}^M \left( H(X^*_m)^{t_m}   \,S^{ n_m}_{t_m,a} \right)~,
\eeqn
where $S^N_{M,a}$ is defined in \req{defS}.
\end{corollary}
Notice if one is Gibbs sampling with the latent
multiplicites $t_m$, then one needs to
sample $t_m$ for all values from 1 up to $n_m$.
This can be a problem if $n_m=1000000$.

\subsection{Table indicators}
A second representation developed in
\cite{ChDuBu11}
stores a table indicator for
each data item in a way that makes it exchangeable.
For this, let the table indicator
$r_n$ be zero if the data $X_n$ does not contribute
a new table, and one if it does contribute a new
table.  The multiplicity $t_k$ for the data value $X^*_k$
is then computed as
the sum for those with the same data value,
so $t_k=\sum_{n=1}^N r_n 1_{X_n=X^*_k} $.
Since we are invariant as to which of the data
starts a table, and there are
$C^{n_{k}}_{t_{k}}$ choices,
the above posterior is modified to yield:
\begin{corollary}[Evidence for discrete case with table indicators]
\label{lem-ev3}
Following the situation of Corollary~\ref{lem-ev2},
\begin{equation}
\label{eq-hpti}
p(X_1,X_2,...,X_N,r_1,...,r_N\,|\,a,b,H(\cdot))
  ~=~\frac{(b|a)_T}{(b)_N}\prod_{m = 1}^M
   \left( H(X_m^{*})^{t_m}S_{t_m, a}^{n_m}
        \frac{1}{{n_{m} \choose t_m}} \right)~,
\end{equation}
where the $t_m$ are derived from the indicators $r_n$.
\end{corollary}
Note the table indicators do not appear explicitly in this other
than through $t_m$ and thus we can ``forget'' the table indicators
and randomly resample them as needed at any stage of Gibbs sampling
since by symmetry their probability of being 1 is $t_m/n_m$.

This representation has two major advantages: one only
needs to incrementally change the $t_m$, not explore all possible values
between $1$ and $n_m$,
and it only requires the use of ratios of Stirling numbers
(for instance $S_{t_m+1, a}^{n_m}/S_{t_m, a}^{n_m}$)
which can be easier to compute.

\subsection{Moments}
Useful quantities to understand the application of the PDP to a
discrete base distribution,
especially for the hierarchical case, are its moments.
We give them here so we can properly interpret the discrete case.

\begin{lemma}[Moments for the discrete case]\label{lem-dcm}
Assume the discrete base distribution $H(\cdot)$ is over the
integers $\SetN$, with probability vector $\vec\theta$, so there is
probability $\theta_k$ for the value $k$. Let $\vec{p}~\sim~
\mbox{PDP}(a,b,H)$. Then the mean, variance, covariance and third
order moments of $\vec{p}$ according to this prior are given by
\bqan
  \expec{}{\vec{p}} &=& \vec{\theta} ~.\\
  \vari{}{p_k} &=&  \frac{1-a}{b+1} \theta_k (1- \theta_k) \\
  \covi{}{p_{k_1}}{p_{k_2}} &=&  - \frac{1-a}{b+1} \theta_{k_1}  \theta_{k_2} \vspace*{-6pt}
\eqan
\bqan
  \lefteqn{
  \expec{}{(p_{k_1}-\theta_{k_1})(p_{k_2}-\theta_{k_2})(p_{k_3}-\theta_{k_2})} }\\
  &=& \left\{
  \begin{array}{lr}
  2\frac{(1-a)(2-a)}{(b+1)(b+2)}\theta_{k_1}\theta_{k_2}\theta_{k_3} &
  \mbox{when $k_1,k_2,k_3$ disjoint} \\
  \frac{(1-a)(2-a)}{(b+1)(b+2)}(2\theta_{k_1}-1)\theta_{k_1}\theta_{k_2} &
  \mbox{when } k_1=k_2\neq k_3 \\
  \frac{(1-a)(2-a)}{(b+1)(b+2)} \theta_{k_1}(1-\theta_{k_1})(1-2\theta_{k_1}) &
  \mbox{when } k_1=k_2=k_3
  \end{array}
  \right.
\eqan
\end{lemma}
Now consider the case where $H(\cdot)$ has domain $1,...,K$, and
probability vector $\vec{\theta}$. Denote this by
$\mbox{discrete}(\vec\theta)$. Consider a $K$ dimensional Dirichlet
distribution with parameters given by $\alpha\vec{\theta}$. This has
corresponding moments
\bqan
  \expec{}{\vec{p}} &=& \vec{\theta} ~.\\
  \vari{}{p_k} &=&  \frac{1}{\alpha+1} \theta_k (1- \theta_k) \\
  \covi{}{p_{k_1}}{p_{k_2}} &=&  - \frac{1}{\alpha+1} \theta_{k_1}  \theta_{k_2} \vspace*{-6pt}
\eqan
\bqan
  \lefteqn{
  \expec{}{(p_{k_1}-\theta_{k_1})(p_{k_2}-\theta_{k_2})(p_{k_3}-\theta_{k_3})}} &&\\
  &=& \left\{
  \begin{array}{lr}
  \frac{4}{(\alpha+1)(\alpha+2)}\theta_{k_1}\theta_{k_2}\theta_{k_3} &
  \mbox{when $k_1,k_2,k_3$ disjoint} \\
  \frac{2}{(\alpha+1)(\alpha+2)}(2\theta_{k_1}-1)\theta_{k_1}\theta_{k_2} &
  \mbox{when } k_1=k_2\neq k_3 \\
  \frac{2}{(\alpha+1)(\alpha+2)} \theta_{k_1}(1-\theta_{k_1})(1-2\theta_{k_1}) &
  \mbox{when } k_1=k_2=k_3
  \end{array}
  \right.
\eqan
Thus we can conclude following:   When $0<a\ll 1$, then we have that
$\mbox{PDP}(a,b,\mbox{Categorical}(\vec\theta))$ is
approximated by a $\mbox{Dirichlet}\left(\frac{a+b}{1-a}\vec\theta\right)$
(and it is already known equality holds when $a=0$).
The two
distributions differ by a factor of $O(a^2)$ in all the
moments of order one to three.
Thus the PDP applied to finite discrete
distributions is approximated by a proper Dirichlet.

\subsection{Computing Stirling numbers}
To work with the discrete case, one needs to sample the
multiplicities $t_m$.
These can be sampled using Gibbs sampling and precomputed tables of
the Stirling numbers $S^{ n}_{t,a}$.
When used in sampling, the Stirling numbers $S_{t,a}^n$
need to be tabulated or cached for the required discount parameter $a$.
Because they rapidly become very large, they need to be stored
and computed in log format.
The recursion must be used rather than the exact formulation,
and becomes
\newcommand{\logS}{\mbox{log}S}
\begin{equation}
\label{logrec}
\logS^{n+1}_{t,a} = \logS^{n}_{t,a}
   + \log \left(  \exp\left(\logS^{n}_{t-1,a}-\logS^{n}_{t,a}\right)
      + (n-t\,a) \right) ~.
\end{equation}
The $\log()$ and $\exp()$ functions can make the evaluation
slow if implemented naively.
Moreover, memory requirements for the full table
of $S_{t,a}^n$ with a fixed discount $a$ is
$n(n+1)/2$ floats for $t\leq n$.
When keeping the discount $a$ fixed, we can reduce memory
with two tricks:
(1)
placing a maximum value on $t$ say 100 or 1000
to limit the cache
and (2)
striping the cache for higher values
of $t$, so only every $L$-th value is stored,
entries of $S_{t, a}^{n}$
for $n=N, N+L, N+2L,...$ and $t<n$.
Computational  time in this case becomes $O(2^L/L)$.
These two techniques save considerable memory at
a small factor in computational time or sampling accuracy.

\subsection{Ratios of Stirling numbers}
When sampling with the table indicator representation, one repeatedly needs
ratios of Stirling numbers.  Denote the Stirling number ratios
\begin{equation}
U_{t,a}^n ~=~\frac{S_{t,a}^{n+1}}{S_{t,a}^n}
~~~~~~~~~~~~~~~~~~~~~~~~
V_{t,a}^n ~=~ \frac{S_{t,a}^{n}}{S_{t-1,a}^n} ~.
\end{equation}
These have the advantage that they do not need to be stored in
log space.  The first ratio, $U_{t,a}^n$,
is readily computed from the second,  $V_{t,a}^n$,
so is not stored:
\begin{eqnarray*}
U_{1,a}^n &=& n-a~~~~~~~~~~~~~~~~~~~~~~~~~~\mbox{ for $n\ge 1$}\\
U_{t,a}^n &=& \frac{1}{V_{t,a}^n} + ( n-t\,a)~~~~~~~\mbox{ for $n\ge t >1$} ~.
\end{eqnarray*}
The second formula follows directly from the linear recursion
for $S_{t,a}^{n}$.
The second ratio, $V_{t,a}^n$, has the following recursion
\begin{eqnarray}
\label{ratiorec}
V_{n,a}^{n} &=& \frac{1}{U_{n-1,a}^{n-1}}
~~~~~~~~~~~~~~~~~~~~~~~~~~~~~~~~~~~~\mbox{ for $n \ge 2$} \\
\nonumber
V_{t,a}^{n+1} &=& \frac{1}{U_{t-1,a}^n}\left(1 + ( n-t\,a)V_{t,a}^{n}\right)
~~~~~~~\mbox{ for $n\ge t \ge 2$} ~.
\end{eqnarray}
Thus a simple recursion yields tables for
$V_{t,a}^{n}$ from which $U_{t,a}^n$ can be computed
without any resort to the log storage of
the Stirling numbers or the use of transcendental functions.

Moreover, these ratios yield more
accurate computation because the direct log computation of the
Stirling numbers leads to some loss of precision.
Values in Table~\ref{tbl:cmp} computed with
discount parameter $a=0.5$ presents the results using floating point
(32 bit)
calculation versus the double precision results (64 bit)
as a proxy for evaluating round-off error.
\begin{table}
\begin{center}
\caption{Comparative values for ratio versus log calculations for $a=0.5$}
\label{tbl:cmp}
\begin{tabular}{|c|c|c|c|}
\hline
$(n,t)$ & double $V^n_{t,a}$ & float $V^n_{t,a}$ & float $\exp\left(\logS^{n}_{t,a}-\logS^{n}_{t-1,a}\right)$
 \\
\hline
(10000,10) & 0.222133 & 0.222124 & 0.180696 \\
(10000,100) &0.0201025 & 0.0201025 &  0.0262359 \\
(10000,1000) & 0.00189684 & 0.00189685 & 0.00181349 \\
\hline
\end{tabular}
\end{center}
\end{table}
One can see that the linear recursion in log space,
Equation~(\ref{logrec}),
yields considerably less accurate results than
the ratio recursion of Equation~(\ref{ratiorec})  when
computation is done with floats rather than doubles.
Computation using the ratios can also be upto 5 times faster.

\section{Discussion}\label{sct-conc}

For the non-atomic case of the two parameter Poisson-Dirichlet
distribution, consistency, convergence and posterior results have
been presented, mostly drawn from the literature, though some proofs
are given in the Appendix. We have augmented these results with a
number of plots to illustrate the nature of the underlying
distributions.   Most significantly, in
practice we recommend fitting one of the
two parameters $a$ or $b$ of the PDP or PDD when performing inference
with them.

\paragraph{Chinese Restaurant distribution:}
The distribution on partitions induced by the two parameter Poisson-Dirichlet
distribution, called here the Chinese Restaurant distribution,
has also been presented, together with a summary
of its use in a scheme for generating trees.
The fragmentation and coagulation
duality for the Chinese Restaurant distribution and the
two parameter Poisson-Dirichlet
distribution followed: coagulation allows simplification
of some hierarchical Poisson-Dirichlet distributions, and
fragmentation allows simplification
of some nested Poisson-Dirichlet distributions.

\paragraph{Improper Priors:}
The infinite probability vector underlying the Poisson-Dirichlet
distribution was shown to be in the form of an improper
Dirichlet on any finite sub-vector.  As soon as data is presented,
the posterior becomes proper (or normalised).
Alternatively, if a size-biased ordering of the probabilities in the
infinite vector is made, the knowledge implicit in
making a size-biased order turns the improper Dirichlet
into a GEM distribution, the standard
``stick-breaking'' distribution most commonly used to define a
Poisson-Dirichlet distribution.

\paragraph{Discrete distributions:}
For the discrete case, not well covered in the Probability and
Statistics literature, posterior results have also been presented.
Moreover, it has been shown that the two parameter Poisson-Dirichlet
Process with discount $a>0$ and concentration $b$
on a discrete base distribution behaves rather like a
Dirichlet distribution with concentration
$\frac{a+b}{1+a}$.
The Dirichlet Process on a discrete base distribution is identical to
a Dirichlet distribution.

\paragraph{Semi-conjugacy:}
While the posteriors presented for the
discrete case of the PDP introduced latent values,
multiplicities or indicators, they are
conjugate to the multinomial family.  This means
the Poisson-Dirichlet Process and the Dirichlet Process
can be used hierarchically and
remarkably a form of semi-conjugacy applies
({\it i.e.}, conjugacy at the expense of introducing latent variables).

\paragraph{Second order Stirling numbers:}
Computation of the Stirling numbers needed to perform Gibbs
sampling was also presented.  The closed form for
computing Stirling numbers looks like a difference approximation
to an $M$-th order derivative, thus it is intrinsically unstable to
compute for higher values of $M$.  The standard linear recursion
therefore seems necessary, although double precision
is needed for larger values.  Ratios of Stirling numbers, however,
can be more accurately computed.

\paragraph{Acknowledgements.}
NICTA is funded by the Australian Government as represented by the
Department of Broadband, Communications and the Digital Economy and
the Australian Research Council through the ICT Centre of Excellence
program.


\bibliographystyle{alpha}
\newcommand{\etalchar}[1]{$^{#1}$}

\newpage
\appendix
\section{Proofs}
\subsection{Proofs for Section~\ref{sct-pp}}

\subsubsection{Proof of Lemma~\ref{lem-exp}}

Build on the result from Lemma~\ref{lem-ev} using
Definition~\ref{defn-pis}. The formula of Lemma~\ref{lem-ev} also
applies to $I^*_N$, the indices $I_N$ with size-biased ordering
applied, so $p(I^*_N)$, but remove the terms in
$H(X^\#_k)$. Using the notation of Definition~\ref{defn-pis} and
this lemma, we get the form
\beq\label{eqn-l3}
  p(I^*_N)~= ~ \frac{(b|a)_M}{(b)_N}
   \prod_{m=1}^M \frac{\Gamma(n_m-a)  }{\Gamma(1-a)}~.
\eeq
Now one can marginalise out the entires in $I^*_N$ but keeping the
constraint that there are $M$ distinct $k$'s in there, which will
affect the last product of $M$ terms only.

The indexes $1,...,M$ occur in the sequence $I^*_N$ of size $N$.
Ignoring the ordering constraints of size-biased ordering, there
are $N$ choose $n_1,\, ..., n_M$, $ C^{N}_{n_1,\, ..., n_M} $ ways
the indexes can occur in $I^*_N$. Now adjust this for the ordering
constraints. For every sequence starting with $1$ there exists some
starting with $2$, ..., $M$. By symmetry, $n_1/N$ of the sequences
start with $1$. Now how many of these have the second integer
appearing in sequence being $2$?   Again by symmetry, $n_2/(N-n_1)$
of the sequences starting with $1$ have $2$ as the next integer in
sequence. Likewise, of those sequences with $1$, $2$ being the first
two occurring integers respectively, $n_3/(N-n_1-n_2)$ have $3$
occurring next. Thus, the number of sequences by size-biased ordering
with counts $n_1,\, ..., n_M$ are
\beqn
  C^{N}_{n_1,\, ..., n_M} \prod_{m=1}^M \frac{n_m}{N-\sum_{i=1}^{m-1} n_i} ~.
\eeqn
Inspection shows this evaluates to an integer since each term
$N-\sum_{i=1}^{m-1} n_i$ divides into $N!$.

To marginalise out the indexes $I^*_N$ in \req{eqn-l3} then, one does
\bqan
  p(M|N) &=&
  \sum_{\sum_{m=1}^M n_m = N,\,\,n_m\ge 1}
  p(I^*_N)~C^{N}_{n_1,\, ..., n_M} \prod_{m=1}^M \frac{n_m}{N-\sum_{i=1}^{m-1} n_i} ~.
\\
  &=& \frac{(b|a)_M}{(b)_N}
  \quad N!\nq\nq \sum_{\sum_{m=1}^M n_m = N,~ n_m\ge 1}
            \prod_{m=1}^M \left(\frac{\Gamma(n_m-a)  }{\Gamma(n_m+1)\Gamma(1-a)}
                 \frac{n_m}{N-\sum_{i=1}^{m-1} n_i} \right)~.
\eqan
The full summation formula for $S^N_{M,a}$ follows. \qed

\subsubsection{Generalized Stirling Numbers}\label{appGS}

We need the following expressions for generalized Stirling numbers.
All but the explicit expression (iii) are due to  \cite{Hsu:88}.

\begin{theorem}[Expressions for Generalized Stirling Numbers]\label{thmGS}
The following expressions all define the same generalized Stirling
numbers $S(n,k;\a,\b,r)$, where the parameters $\a,\b,r\in\SetR$
have been suppressed when constant.

\begin{tabbing}
\phantom{(ii)} \= \phantom{Explicit expressions} \= \kill
(o) \> Implicit: \>
$\displaystyle(t|-\a)_n = \sum_{k=0}^n S(n,k)(t-r|-\b)_k$\\
\> Both sides are polynomials in $t$ of degree $n$. $(z|a)_n:=z(z\!+\!a)...(z\!+\!(n\!-\!1)a)$.
\\
(i) \> Linear recursion: \>
$\displaystyle S(n\!+\!1,k) = S(n,k\!-\!1) + (k\b-n\a+r) S(n,k)$
\\
    \> Boundary cond.: \>
$ S(n,k)=0$ for $k>n$, $\quad S(n,0)=(r|-\a)_n$
\\
(ii) \> Mult.\ recursion: \>
${N\choose K} S(N,K,\a,\b,R) = $ \hspace{25ex} (any $k,r$)\\
\> \> $= {\displaystyle\sum_{n=0}^N} {N\choose n} S(n,k;\a,\b,r) S(N\!-\!n,K\!-\!k;\a,\b,R-r)$
\\
(iii) \> Explicit expression: \>
$\displaystyle S(n,k) = {1\over k!\,\b^k}\sum_{j=0}^k{\textstyle{k\choose j}}(-)^{k-j}(\b j+r|-\a)_n\quad\qquad$ ($\b\neq 0$)
\\
(iv) \> Generative fct.: \>
$\displaystyle\sum_{n=0}^\infty S(n,k){t^n\over n!} = {(1+\a t)^{r/\a}\over k!}\left({(1+\a t)^{\b/\a}-1\over\b}\right)^k\;\;$ ($\a\b\neq 0$)
\end{tabbing}
\end{theorem}

\paragraph{Proof.}
The generalized Stirling numbers are defined in \cite{Hsu:88} by
(o). Hsu and Liu derive expressions (i),(ii), and (iv) from (o):
It is easy to verify that recursion (i) satisfies definition (o).
Using (i), one can see that the generating function
$g_k(t):=\sum_{n=0}^\infty S(n,k)t^n/n!$ satisfies the differential
equation system
\beqn
  (1+\a t){d\over dt}g_k(t)=g_{k-1}(t)+(k\b+r)g_k(t) \qmbox{with}
  g_k(0)=0 \qmbox{and} g_0(t)=(1+\a t)^{r/\a},
\eeqn
which has a unique solution. Substituting (iv) into this dgl shows
that (iv) is a solution.
If we take a product of the generating functions of
$ S(n,k;\a,\b,r)$ and $ S(N-n,K-k;\a,\b,R-r)$, use (iv), and
identify the coefficients of $t^n$ in both sides, we arrive at the
multiplicative recursion (ii).

Interestingly, Hsu and Liu do {\em not} derive the explicit expression
(iii), although it easily follows by Taylor expanding the r.h.s.\ of
(iv) and by identifying the coefficients of $t^n$ as follows:
The binomial identity gives
\beqn
  ((1\!+\!\a t)^{\b/\a}-1)^k \;=\; \sum_{j=0}^k{\textstyle{k\choose j}} (-)^{k-j}(1\!+\!\a t)^{\b j/\a}
\eeqn
Exploiting this, (iv), and
$(1+z)^\gamma=\sum_{n=0}^\infty{\gamma\choose n}z^n$, where
${\gamma\choose n}={\Gamma(\gamma+1)\over n!\,\Gamma(\gamma-n+1)}$, we
get
\bqan
  \sum_{n=0}^\infty S(n,k){t^n\over n!} &=&
  {1\over k!\,\b^k}\sum_{j=0}^k{\textstyle{k\choose j}}(-)^{k-j}(1\!+\!\a t)^{\b j+r\over\a}
\\
  &=& {1\over k!\,\b^k}\sum_{j=0}^k{\textstyle{k\choose j}}(-)^{k-j}\sum_{n=0}^\infty{{\b j+r\over\a}\choose n}(\a t)^n
\\
  &=& \sum_{n=0}^\infty {1\over k!\,\b^k}\sum_{j=0}^k{\textstyle{k\choose j}}(-)^{k-j}(\b j+r|-\a)_n{t^n\over n!}
\eqan
\qed

\subsubsection{Proof of Theorem \ref{thmSN}}\label{appSN}

The proof is based on (a) a recursion for $p(M|N)$, (b) the
expressions for the generalized Stirling numbers in Appendix~\ref{appGS},
and of course (c) the definition \req{defS} of
$S_{M,a}^N$. In order to distinguish between different $M$ as the
sample size $N$ increases, use $M_N$ to denote the value at sample
size $N$.

\paragraph{\bf (i)} 
We exploit recursion
\beqn
  p(M_{N+1}=m|M_N) \;=\; 1\!\!1_{M_N=m-1} {b+(m\!-\!1)a\over b+N} +
                      1\!\!1_{M_N=m} {N-ma\over b+N} ~,
\eeqn
which easily follows from the predictive sampling distribution
\req{eq-smp}.  Multiplying each side by $p(M_N)$, and summing
over $M_N$ this becomes
\beqn
  p(M_{N+1}=m) \;=\; p(M_N=m-1){b+(m\!-\!1)a\over b+N} +
                      p(M_N=m){N-ma\over b+N}
\eeqn
Inserting the explicit expression
$p(M_N=m)= S_{m,a}^N(b|a)_m/(b)_N$
of Lemma \ref{lem-exp} into this recursion and canceling all common
factors we get
\beqn
  S_{m,a}^{N+1} \;=\; S_{m-1,a}^N + (N-ma) S_{m,a}^N.
\eeqn
The boundary conditions $ S_{m,a}^N=0$ for
$m>N$ and $ S_{0,a}^N=\delta_{N,0}$ follow from the explicit
expression in Definition \req{defS} or simply by reflecting on the
meaning of $p(M_N=m)$.

\paragraph{\bf (ii) and (iii)} 
Consider the generalized Stirling numbers $S(n,k;\a,\b,r)$ for the
special choice of parameters $(\a,\b,r)=(-1,-a,0)$. For this choice,
recursion (i) of Theorem \ref{thmGS} reduces to recursion (i) of
Theorem \ref{thmSN}, including the boundary conditions. Hence
$S_{M,a}^N=S(N,M;-1,-a,0)$.

It is easy to see that also (ii) and (iii) of Theorem \ref{thmGS}
reduce to the first expression of (ii) and (iii) of Theorem
\ref{thmSN} for $(\a,\b,r)=(-1,-a,0)$, which shows that the
expressions are equivalent.

The special case for (iii) when $a=0$ holds by noting
that (iii) for $a>0$ is in fact an $M$-point interpolation
to an $M-1$-th partial derivative.  As $a\rightarrow 0$,
this becomes the partial derivative.

The last expression in (ii) follows from the definition of
$S_{M,a}^N$ in \req{defS} by splitting the sum into
$\sum_{n_1=1}^{N-M+1}$ and $\sum_{n_2+...+n_M=N-n_1}$ and the
product into $m=1$ and $m>1$, and identifying the terms with
${N-1\choose n_1-1}$, $S_{1,a}^{n_1}$ and $S_{M-1,a}^{N-{n_1}}$.

Note that the first expression in (ii) does {\em not} reduce to the
second expression for $m=1$. Nevertheless, the (very different!)
derivations of the two expressions show that they must be equal.

\paragraph{\bf (iv)} 
Using $\Gamma(N+x)/\Gamma(N+y)\simeq N^{x-y} $ for large $N$,
we see that the $m$-contribution in (iii)
is asymptotically proportional to
\beqn
  \prod_{h=0}^{N-1}(h-am) \;=\; {\Gamma(N-am)\over\Gamma(-am)}
  \;=\; {-am\Gamma(N)\over\Gamma(1-am)}{\Gamma(N-am)\over\Gamma(N)}
  \;\;\stackrel{N\to\infty}{\simeq}\;\; {-am\Gamma(N)\over\Gamma(1-am)}{1\over N^{am}}
\eeqn
Due to the factor $m$, the $m=0$ term does not contribute.
So the dominant contribution is from $m=1$, followed by $m=2$, etc.
The $m=1$ term yields
\beqn\label{SaAsympt2}
  S_{M,a}^N \;\simeq\; {1\over M! a^M}M{a\Gamma(N)\over\Gamma(1-a)}{1\over N^a}
  \;=\; {1\over\Gamma(1\!-\!a)}{1\over\Gamma(M)\,a^{M-1}}{\Gamma(N)\over N^a}
\eeqn
The relative accuracy is $O(M/N^a)$, i.e.\ the approximation is good
for $M\ll N^a$. The smaller $a$, the larger $N$ needs to be to reach
a reasonable accuracy. Higher $m$-terms may be added, but the
alternating sign indicates cancelations and hence potential
numerical problems.

\paragraph{\bf (v)} 
follows from $S_{M,0}^N=S(n,k;-1,0,0)=|S(n,k;1,0,0)|$ and the fact
that $S(n,k;1,0,0)$ are Stirling numbers of the first kind from
\cite{Hsu:88}. \qed

\subsubsection{Proof of Lemma~\ref{lem-part}}

We need to differentiate the different $M$ that results from the
partition sample $I^*_N$ as $N$ increases. Subscript $M$ as $M_N$ so
we can differentiate it as $N$ changes. When $M_N$ is known, the
following series relation holds:
\beqn
  \expec{M_N}{M_{N+1}} ~=~ \frac{b+M_N a}{N+b} + M_N
  ~=~ \frac{b}{N+b} + \frac{a+ b+N}{N+b} M_N ~.
\eeqn
Taking expected values across $M_N$ yields the recursive form
\beqn
  \expec{}{M_{N+1}} ~=~ \frac{b}{N+b} + \frac{a+ b+N}{N+b} \expec{}{M_N}
\eeqn
The equation for $\expec{}{M_N} $ given in the lemma is proven from
this by induction, with the value $1$ when $N=1$. Note the
derivation of the solution to the above recursive formula was made
by unfolding the recursion into a summation, and then simplifying
the summation using hypergeometric functions.

The approximation for $\expec{}{M_N} $ given in the lemma is derived
for $N,b\gg a$ as follows:
\bqan
\frac{(a+b)_N}{(b)_N} &=& \exp\left(\log\Gamma(a+b+N)-\log\Gamma(b+N)
       -\left( \log\Gamma(a+b)-\log\Gamma(b) \right)\right) \\
  &\simeq& \exp\left( a \left( \psi_0(b+N)-\psi_0(b) \right) \right)\\
  &\simeq& \left(1+\frac{N}{b}\right)^a \exp\left( \frac{-a}{2} \left(\frac{1}{b+N}- \frac{1}{b}\right)  \right) \\
  &=&  \left(1+\frac{N}{b}\right)^a \exp\left( \frac{aN}{2b(b+N)}   \right) ~.
\eqan
The first approximation step makes a first order Taylor expansion
since $a$ is small, $0<a<1$, and the second approximation step uses
an approximation for $\psi_0(b) $ with error $O(1/b^2)$.

For the expected variance, a similar strategy is used but the steps
are  more complicated. The following series relation holds:
\bqan
\expec{M_N}{M_{N+1}^2} &=&
  \frac{b+M_N a}{N+b}(M_N+1)^2  + \frac{N-M_N a}{N+b}M_N^2 \\
 &=&  \frac{b}{N+b} + \frac{2b+a}{N+b} M_N
  + \frac{2a+ b+N}{N+b} M^2_N ~,
\eqan
where $\expec{}{M_N^2}=1$ when $N=1$.
Taking expected values over $M_N$ yields the recursive form
\beqn
  \expec{}{M_{N+1}^2} ~=~ \frac{b}{N+b} + \frac{2b+a}{N+b}
  \expec{}{M_N}
    + \frac{2a+ b+N}{N+b} \expec{}{M_N^2} ~.
\eeqn
Evaluation of this recursive formula can be made as before,
and the result is the formula
\beqn
  \expec{}{M_N^2} ~=~ \frac{b(a+b)}{a^2}\frac{(2a+b)_N}{(b)_N}
   - \frac{b(2b+a)}{a^2}\frac{(a+b)_N}{(b)_N} + \frac{b^2}{a^2} ~.
\eeqn
The result then comes from evaluating $\expec{}{M_N^2} -
\left(\expec{}{M_N}\right)^2$ and simplifying terms. The
approximation proceeds as before.

To handle the case where $a=0$ the same recursive formula for
$\expec{}{M_N}$ and $\expec{}{M_N^2}$ can be used, but are evaluated
differently since $a=0$.  The closed form formula for
$\expec{}{M_N}$ follows clearly by induction on $N$.  The closed
form formula for  $\expec{}{M_N^2}$, readily proven by induction, is
\beqn
  b(\psi_0(b+N)-\psi_0(b)) +b^2\left(\psi_0(b+N)-\psi_0(b)\right)^2
  + b^2(\psi_1(b+N)-\psi_1(b))~.
\eeqn
Subtracting off $\left(\expec{}{M_N}\right)^2$ yields the result for
$\vari{}{M_N}$. The approximations for both $\expec{}{M_N}$ and
$\vari{}{M_N}$ follow by taking the first order terms of
$\psi_0(\cdot)$ and $\psi_1(\cdot)$, the log and the inverse
respectively. \qed

\subsubsection{Proof of Lemma~\ref{lem-series}}

The value $M$ is equal to the number of indices that have a non-zero
count in the sample of size $N$. Given probability vector $\vec q$,
the probability that index $k$ has a non-zero count after $N$
samples is $1-(1-q_k)^N$. Summing these over all $k$ gives an upper
bound.

To generate  bounds on $1-(1-q_k)^N$, note $1-(1-q_k)^N\leq 1$, and
the bound is closer the larger $q_k$. Second, by Taylor expansion
\beqn
  1-(1-q_k)^N ~=~ N q_k - \frac{N(N-1)}{2} q_k'^2
\eeqn
for some $0\leq q_k' \leq q_k $. So $1-(1-q_k)^N \leq N q_k$, and
the bound is closer the smaller $q_k$, especially when $N q_k\ll 1$.
Put these two bounds together and we get for any positive integer
$m$
\beqn
  \sum_{k=1}^\infty 1-(1-q_k)^N ~\leq~
    m + N \sum_{k=m+1}^\infty q_k~.
\eeqn

For the geometric series, $q_k=r^{k-1}(1-r)$, the sum in the bound
evaluates to $r^m$, so we seek to minimise $m+Nr^m$. This bound can
be modified if we let $m\in\SetR^+$ to
\beqn
  \min_{m\in \SetN^+} m+Nr^m ~\leq~
  \min_{m\in\SetR^+} 1+m+Nr^{m-1}
\eeqn
Differentiating yields a minima at $r^{m-1}=\frac{1}{N \log 1/r}$.
The result follows by substitution.

For the Dirichlet series, $q_k=\frac{k^{-s}}{\zeta(s)}$,
the sum in the bound can be bounded by an integral
\bqan
  \sum_{k=m+1}^\infty q_k &\leq& \int_{m+1/2}^\infty
      \frac{1}{(k)^{s}\zeta(s)}~, \\
  &=& \left. \frac{-1}{(k)^{s-1}\zeta(s)(s-1)} \right|_{m+1/2}^\infty \\
  &=& \frac{1}{(m+1/2)^{s-1}\zeta(s)(s-1)}
\eqan
As before, modifying the bounds yields the formula to minimise
\beqn
  m+1+ \frac{N}{(m-1/2)^{s-1}\zeta(s)(s-1)} ~.
\eeqn
Differentiation gives a minimum at $N = (m-1/2)^{s}\zeta(s)$ and so
the bound follows. \qed

\subsection{Proofs for Section~\ref{sct-fc}}

\subsubsection{Proof of Theorem~\ref{th-pfrag}}

One sees the samples generated by the fragmentation process.
Given a partition with partition count $L$
with occurrences $n_1,...,n_L$ totalling $N$,
one needs to assign every entry $l$ to a latent
cluster $k_l$.  This results in a coarser partition
with count $K$ bins, and with total occurrence $m_k=\sum_{l\,:\,k_l=k}n_l$
(where  $m_k=|P_m|$ for the latent bin $P_m$)
from $t_k=\sum_{l\,:\,k_l=k}1$ original bins,
totalling $\sum_{k=1}^K t_k=L$.
The full probability for the fragmentation components is
(note $a_2$ could also be zero here, but these terms correctly cancel,
so the parallel equations for $a_2=0$ are not given):
\bqan
&&
\frac{(b|a_1a_2)_K}{(b)_N}
\prod_{k=1}^K(1-a_1a_2)_{m_k-1} ~\cdot~
\prod_{k=1}^K \frac{(-a_1a_2|a_1)_{t_k}}{(-a_1a_2)_{m_k}}
\prod_{l\,:\,k_l=k}(1-a_1)_{n_l-1}\\
&=& \frac{1}{(b)_N}\prod_{l=1}^L(1-a_1)_{n_l-1}
(b|a_1a_2)_K
\prod_{k=1}^K \frac{(-a_1a_2|a_1)_{t_k}}{-a_1a_2}\\
&=& \frac{1}{(b)_N} a_1^L \prod_{l=1}^L(1-a_1)_{n_l-1}
\cdot
(b/a_1|a_2)_K\prod_{k=1}^K (1-a_2)_{t_k-1} ~.
\eqan
We need to marginalise over all possible coarser partitions
or assignments to $k_l$.    Now the second term is in the form
of a $\mbox{CRD}\left(\{1,...,l\},a_2,b/a_1\right)$
so marginalising out the $t_k$ yields $(b/a_1)_L$,
which multiplied by $a_1^L$ gives $(b|a_1)_L$.
The result is that the above probability becomes
\[
 \frac{(b|a_1)_L}{(b)_N}  \prod_{l=1}^L(1-a_1)_{n_l-1} ~.
\]
This is the $\mbox{CRD}\left(\{1,...,N\},a_1,b\right)$.

\subsubsection{Proof of Theorem~\ref{th-pcoag}}

Consider a sample partition with a bin $j$ with $n_j$ entries,
however we do not know which sub-partitions
(from the unknown index $k$) the entry comes from,
nor how many there were.  So now each bin $j$ would be sub-partitioned
into $L_j$ bins with occurrence counts $m_{j,l}$ each where
$\sum_{l=1}^{L_j}m_{j,l}=n_j$. In total, there are $L=\sum_{j=1}^JL_j$ bins.
With this assignment done, the probability for the coagulation
becomes:
\begin{eqnarray*}
&& \frac{(b|a_1)_L}{(b)_N} \prod_{j,l} (1-a_1)_{m_{j,l}-1} \cdot
 \frac{(b/a_1|a_2)_J}{(b/a_1)_L} \prod_{j} (1-a_2)_{L_j-1} \\
&=& \frac{(b|a_1a_2)_J}{(b)_N}  \frac{1}{(-a_1a_2)^J}
   \prod_{j=1}^J    (-a_1a_2|a_1)_{L_j} \prod_{l=1}^{L_j} (1-a_1)_{m_{j,l}-1} ~.
\end{eqnarray*}
Note the last formula is undefined for $a_2=0$, however, the
zero terms cancel and a parallel proof for the $a_2=0$ case
is readily seen to hold.
To construct the marginalised probability, we have to
marginalise over all possible sub-partitions
($L_j$ bins with occurrence counts $m_{j,l}$).
The term inside the $\prod_{j=1}^J$ is the form of a
$\mbox{CRD}\left(\{1,...,L_j\},a_1,-a_1a_2\right)$ so marginalising out the
partitions (represented by $L_j$ and $m_{j,l}$)
yields $(-a_1a_2)_{n_j}=(1-a_1a_2)_{n_j-1}(-a_1a_2)$
(again, noting a similar argument applies for the $a_2=0$ case).
The probability simplifies then to
\[
 \frac{(b|a_1a_2)_J}{(b)_N}  \prod_{j=1}^J    (1-a_1a_2)_{n_j-1} ~.
\]
This is then $\mbox{CRD}\left(\{1,...,N\},a_1a_2,b\right)$
as required.

\subsubsection{Proof of Theorem~\ref{th-frag}}

To prove the result, it is sufficient to show for all samples,
the two have equivalent marginalised posteriors.
Consider a sample partition, for
a given bin each entry would be drawn with a probability
$p_kq_{j,k}$, however we do not know the $j$ or $k$ associated.
We have to assign the $k$, whereas the $j$ term is simply a
sub-partition so we do not need to know it.
Given a sample partition like this of $L$ bins
with sizes $n_1,...,n_L$ totalling $N$,
one needs to assign every entry $l$ to a
cluster $k_l$.  This results in a coarser partition
with $K$ bins, with total count $m_k=\sum_{l\,:\,k_l=k}n_l$
from $t_k=\sum_{l\,:\,k_l=k}1$ original bins,
totalling $\sum_{k=1}^K t_k=L$.
But with the clusters $k_l$ assigned, the posterior
with the marginalised terms for $\vec{p}$ and $\vec{q}_k$'s can be
given.
This posterior corresponds to the situation in
Theorem~\ref{th-pfrag}, so apply that result and
this converts the posterior into a
probability for $\mbox{CRD}\left(\{1,...,N\},a_1,b\right)$.
This is the required posterior for the sample from
$\mbox{PDD}\left( a_1, b\right)$.
\qed

\subsubsection{Proof of Theorem~\ref{th-coag}}

As before, to prove the result, it is sufficient to show for all samples,
the two have equivalent marginalised posteriors.
Consider a sample partition with a bin $j$ with $n_j$ entries,
each entry would be drawn with a probability
$\sum_{k\,:\,j_k=j}p_k$, however we do not know which sub-partitions
(from the unknown index $k$) the entry comes from,
nor how many there were.  So now each bin $j$ would be sub-partitioned
into $L_j$ bins with occurrence counts $m_{j,l}$ each where
$\sum_{l=1}^{L_j}m_{j,l}=n_j$. In total, there are $L=\sum_{j=1}^JL_j$ bins.
With this assignment done, the marginalised posterior can be written down
from the term for $\vec{p}$ and the term for $\vec{q}$,
however this corresponds to the situation
from Theorem~\ref{th-pcoag}.
Applying that result yields a marginalised posterior in
the form of a $\mbox{CRD}\left(\{1,...,N\},a_1a_2, b\right)$.
This is the required posterior for the sample from
$\mbox{PDD}\left( a_1a_2, b\right)$.

\subsection{Proofs for Section~\ref{sct-dm}}

\subsubsection{Proof of Lemma~\ref{lem-cons}}

Consider  the prior measure for $p_1,...,p_M,p^+_M$. Do a change of
variables to $p_1,...,p_{{M-1}},q_M,p^+_{M-1}$ where $q_M =
p_M/p^+_{M-1}$ and $p^+_{M-1}=p_M+p^+_M$. The Hessian of this change
is $1/p^+_{M-1}$, and the domain goes from the constraint set
$\{p_1\geq 0, ..., p_M\geq 0, p^+_M\geq 0\}$ to $\{p_1\geq 0, ...,
p_{{M-1}}\geq 0, p^+_{M-1}\geq 0, 0\leq q_M \leq 1\}$. The prior
measure can thus be converted to
\beqn
  \left(q_M^{-a-1} \left(1-q_M\right)^{b+Ma-1}\right)
~ \left(p^+_{M-1} \right)^{b+(M-1)\,a-1} \prod_{m=1}^{M-1} p_{m}^{-a-1} ~,
\eeqn
under the new constraint set. Note the prior measure on  sub-vector
$p_1,...,p_{{M-1}}$, as given in Definition~\ref{defn-pyd}, appears
in the second half of this measure.  The initial part involves only
$q_M$, but its constraints are simply $0\leq q_M \leq 1$ which are
independent of the remaining variables. Thus one is left with a
measure on $p_1,...,p_{{M-1}}$. The measure on the sub-vector is now
consistent with Definition~\ref{defn-pyd}. We can repeat this
process recursively to verify consistency for any other sub-vector.
\qed

\subsubsection{Proof of Theorem~\ref{thm-dap}}

Consider Definitions~\ref{defn-pyd} and~\ref{defn-pis}. Define
$G_{\delta}$  in terms of its projection on the finite sub-spaces
$\{p_1,p_2,...,p_M\}$ for all $M$. Let
\beq\label{eq-dapG}
  p(p_1,p_2,...,p_M,p^+_M)
  ~\propto~
  \left(p^+_M\right)^{b+M\,a-1} \prod_{m=1}^M p_{m}^{-a-1}  ~,
\eeq
where $p^+_M=1-\sum_{m=1}^M p_{m}$ and the domain is constrained to
be $p_m>(1-\sum_{i=1}^{m-1} p_i)\delta$ for $m=1,...,M$, and
$p^+_M\geq 0$. Note that by Definition~\ref{defn-pyd},  $b>-a$, and
thus $b+Ma>0$. Exploiting $p_m>\delta$ for $m=1,...,M$, we show
below that the proportionality constant, i.e.\ the integral over the
constrained simplex, is finite. To show  $G_{\delta}$ is proper, we
need to show that the finite priors for each $M$ are proper and that
consistency holds between these priors for different $M$.

The normalization is done as follows. Use the same change of
variables as in the proof of Lemma~\ref{lem-cons}, however now the
domain is different. The constraint set for the initial variables is
\beqn
  C_{p,M} ~=~ \left\{p_1\geq \delta, ..., p_m\geq \left(1-\sum_{i=1}^{m-1} p_i\right)\delta,
  ..., p_M\geq \left(1-\sum_{i=1}^{M-1} p_i\right)\delta, p^+_M\geq 0\right\}~.
\eeqn
By the change of variables this gets mapped to
\beqn
  C_{q,M} ~=~
  \left\{p_1\geq \delta, ...,  p_{M-1}\geq \left(1-\sum_{i=1}^{M-2} p_i\right)\delta,
  p^+_{M-1}\geq 0, \delta \leq q_M \leq 1\right\} ~.
\eeqn
For the purposes of integration, denote the initial and changed
variable sets as $\vec{p}$ and $\vec{q}$ respectively. Thus the
integration works as follows:
\bqan
  Z_{a,b,M,\delta} &:=& \int_{C_{p,M}} \left(p^+_M\right)^{b+M\,a-1} \prod_{m=1}^M p_{m}^{-a-1} \mbox{d}\,\vec{p} \\
  &=& \int_{C_{q,M}} \left(p^+_{M-1}\right)^{b+(M-1)\,a-1} \prod_{m=1}^{M-1} p_{m}^{-a-1}
      \left( 1- q_M\right)^{ b+M\,a-1 } q_M^{-a-1} \mbox{d}\,q
\\
  &=& Z_{a,b,M-1,\delta} \int_{\delta}^1 \left( 1- q\right)^{ b+M\,a-1 }  q^{-a-1}  \mbox{d}\,q
\\
  &=& ... ~=~ \prod_{m=1}^M  \int_{\delta}^1 \left( 1- q\right)^{ b+m\,a-1 }  q^{-a-1}  \mbox{d}\,q~ .
\eqan
Note this is bounded above by bounding the $q^{-a-1}$ terms from
inside the integral with $\delta^{-a-1}$, and extending the
integrals to the range $[0,1]$. This yields the upper bound
$\delta^{-M(a+1)}\prod_{m=1}^M \frac{\Gamma(b+ma)}{\Gamma(b+ma+1)}$.

Now we prove consistency. We need to show that the projection from
the subset $m=1,...,M$ down to some smaller subset $m=1,...,M'<M$ is
consistent. The change of variables above handled the case where
$p(p_1,p_2,...,p_M,p^+_M)$ was projected down to
$p(p_1,p_2,...,p_{M-1},p^+_{M-1})$. Clearly, the projected prior is
equivalent to the direct definition above (see Lemma \ref{lem-cons}
for details). Thus by induction, one can project the prior from the
subset $m=1,...,M$ down to a any smaller subset $m=1,...,M'<M$, and
get the same prior.  By this condition, and Kolmogorov's Consistency
Theorem, it follows that the prior $G_{\delta}$ exists and is proper
for the full Hilbert space of $\vec{p}$.

Now consider the posteriors for a given sample $I_N$. The posterior
for $I_N$ using the improper prior on PDDs is given in
Lemma~\ref{lem-base}. To deal with the proper prior $G_{\delta}$,
the notion of partition size is needed, as given in
Definition~\ref{defn-pis}. Let $M_N$ be the partition size for a
$I_N$, then $p(p_1,p_2,...,p_M,p^+_M|G_{\delta},I_N)$ is
proportional to
\beqn
  \left(p^+_M\right)^{b+M\,a-1} \prod_{m=M_N+1}^M p_{m}^{-a-1} \prod_{m=1}^{M_N} p_{m}^{n_m-a-1} ~,
\eeqn
where the constraints $C_{p,M}$ hold as before. This is the same
form as the posterior Dirichlet distribution (I) given in
Lemma~\ref{lem-base} where the probabilities are further constrained
by  $C_{p,M}$. The normalizing constant can be worked out as above
to be
\beqn
  Z_{a,b,M_N,\delta} ~=~
  \prod_{m=1}^{M_N} B_{1-\delta}(b+ma,n_m-a)
\eeqn
where $B_x(u,v)=\int_0^x t^{u-1}(1-t)^{v-1}dt$ is the incomplete
Beta function defined for $u,v>0$. In our case, $n_m>0$ for all
$1\leq m \leq M_N$ and $a+b>0$, so the Beta function and incomplete
Beta function are well defined. Note the normalizing constant for
the posterior Dirichlet distribution (I) given in
Lemma~\ref{lem-base} is $Z_{a,b,M_N,0}$.

Now consider the $L_{1}$ distance between the two posteriors,
$p(p_1,p_2,...,p_M,p^+_M|G_{\delta},I_N)$ and the posterior
Dirichlet distribution (I) given in Lemma~\ref{lem-base}.  Note
these differ only in domain. Denote them by $P_\d$ and $P_0$
respectively. Using $P_\d\geq P_0$ on $G_\d$, and $P_\d=0$ on
$G_0\setminus G_\d$, and $\int_{G_\d}P_\d d\vec p=1$, we get
\bqan
\lefteqn{ \frac12 d_1(P_\d,P_0) ~:=~ \sup_A|P_\d[A]-P_0[A]|} &&\\
  &=& \frac12\int_{G_0}|P_\d-P_0|d\vec p \\
  &=& \frac12\int_{G_\d} |P_\d-P_0|d\vec p + \frac12\int_{G_0\setminus G_\d} |P_\d-P_0|d\vec p \\
  &=& \frac12\int_{G_\d} P_\d d\vec p - \frac12\int_{G_\d} P_0 d\vec p + \frac12\int_{G_0\setminus G_\d} P_0 d\vec p \\
  &=& 1-\int_{G_\d} P_0 d\vec p
  \;=\; 1-{Z_{a,b,M_N,\delta} \over Z_{a,b,M_N,0}}\int_{G_\d} P_\d d\vec p \\
 &=& 1-{Z_{a,b,M_N,\delta}\over Z_{a,b,M_N,0}} \to 0 \qmbox{for} \delta\to 0
\eqan
This implies convergence in distribution.
\qed

\subsubsection{Proof of Lemma~\ref{lem-base}}

\paragraph{Proof sketch.}
Note for the Proper Posteriors II claim, since $H(\cdot)$ is
non-atomic, each distinct data $X^*_m$ has a corresponding distinct
index $k^*_m$, thus for the purposes of analysis, assume the indices
are given and w.l.o.g.\ they follow size-biased ordering,
so $k^*_m=m$. Thus to prove the Proper
Posteriors I and II claim about posteriors for $\vec p$, multiply
the prior measure for $(p_1,...,p_M)$ of Definition~\ref{defn-pyd}
by the likelihood, which is in terms of the same sub-vector, and the
posterior measure clearly is proportional to the corresponding
posterior Dirichlet in this lemma. The remaining part of the Proper
Posteriors II claim follows from the model family.

To prove the Sampling claim, note that this just takes the expected
value of the posterior in Proper Posteriors II. To prove the Stick
Breaking claim, note this follows directly from the posterior by
standard properties of the Dirichlet.

The size-biased sampling claim is developed sequentially.
First note $p_1$ has the improper prior
$\mbox{Beta}(-a,b+a)$.
The value $p_2/(1-p_1)$ is {\it apriori} independent of $p_1$
and has a the improper prior
$\mbox{Beta}(-a,b+2a)$.
The value $p_3/(1-p_1-p_2)$ is {\it apriori} independent of $p_1$ and $p_2$
and has a the improper prior $\mbox{Beta}(-a,b+3a)$, {\it etc}.
However, because $p_1$ is size-biased
we  know it is the first in the sample, so we add one
to get a $\mbox{Beta}(1-a,b+a)$.
Likewise, we also know $p_2$ appears first in the sample
(after $p_1$),
since it is sized-biased, so again we add one getting
$\mbox{Beta}(1-a,b+2a)$.
Repeating this yields the standard
stick-breaking definition of the GEM distribution.
\qed

\subsection{Proofs for Section~\ref{sct-post}}

\subsubsection{Proof of Corollary~\ref{lem-ev2}}

\paragraph{Proof sketch.}
In the general case where each draw from $H(\cdot)$ is not
necessarily almost surely distinct, the formula of
Lemma~\ref{lem-ev} also applies to $p(X_1,X_2,...,X_N,k_1,...,k_N)$.
Now one can marginalise out the $k_1,...,k_N$, which will affect the
last product of $M$ terms only.

Given the constraints that $t_m$ represents the multiplicity of
$X^*_k$ and $n_k$ represents the total count of $X^*_k$, then all
values for $k_1,...,k_N$ must be included that satisfy the
constraints. Each $n_k$ will be partitioned into $t_k$ different
indices, each occurring at least once, and totaling $n_k$. Thus the
problem of marginalising out the indices $k_1,...,k_N$ to the
multiplicities $t_1,...,t_M$ is equivalent to the summation over
configurations by size-biased ordering, done for
Lemma~\ref{lem-exp}, and an identical result can be applied. \qed

\subsubsection{Proof of Lemma~\ref{lem-dcm}}

Let $\vec{p}~\sim~ \mbox{PDP}(a,b,H)$. Let $\vec{q}~\sim~
\mbox{PDD}(a,b)$ be the underlying PDD, and let the corresponding
independent samples from $H(\cdot)$ be $X_l \in \cal N$.  From the
definition of a PDP,
\beqn
  p_k ~=~ \sum_{l} q_l 1_{X_l=k}
\eeqn
Taking the expected value of this over $\vec{X}$, yields $ \sum_{l}
q_l \theta_k$, and hence $\theta_k$ irrespective of $\vec{q}$.

Now consider any moment.  We present one case, and others can be
treated similarly. For $k_1, k_2, k_3$ three indices
\bqan
& & \expec{\vec{q},\vec{X}}{\left(  \sum_{l} q_l 1_{X_l=k_1} - \theta_{k_1}  \right)\left(\sum_{l} q_l 1_{X_l=k_2} -  \theta_{k_2}  \right)\left( \sum_{l} q_l 1_{X_l=k_3} - \theta_{k_3} \right)} \\
&=& \expec{\vec{q},\vec{X}}{\sum_{l_1,l_2,l_3} q_{l_1} q_{l_2} q_{l_3} \left(  1_{X_{l_1}=k_1} - \theta_{k_1}  \right)\left( 1_{X_{l_2}=k_2} -  \theta_{k_2} \right)\left(  1_{X_{l_3}=k_3} - \theta_{k_3} \right)} \\
\eqan
Now $X_{l_1}$ is independent of $X_{l_2}$ whenever $l_1\neq l_2$. So
we have to express the sum $\sum_{l_1,l_2,l_3}$ into different equal
and unequal parts so that the expected value over $\vec{X}$ can be
applied. This would be
\beqn
  \sum_{l_1,l_2,l_3} \cdot ~=~ \sum_{l_1,l_2,l_3 \mbox{~disjoint}} \cdot
  + \sum_{l_1=l_2\neq l_3}  \cdot +  \sum_{l_1=l_3\neq l_2}  \cdot +  \sum_{l_2=l_3\neq l_1}  \cdot
  + \sum_{l_1=l_2 = l_3} \cdot
\eeqn
Any sum which has a term with one index, $l_1$ say, not equal to the
others, will contain the expression
\beqn
  \expec{X_{l_1}}{ 1_{X_{l_1}=k_1} - \theta_k } ~=~  \theta_{k_1} - \theta_{k_1} ~=~0~,
\eeqn
and hence can be discarded. Thus for the first three central
moments, the expansion of sums that remains non-zero are
\bqan
  \sum_{l_1,l_2} \cdot  &=& \sum_{l_1=l_2} \cdot \\
  \sum_{l_1,l_2,l_3} \cdot &=& \sum_{l_1=l_2=l_3} \cdot \\
  \sum_{l_1,l_2,l_3,l_4} \cdot &=& \sum_{l_1=l_2\neq l_3=l_4} \cdot  +   \sum_{l_1=l_3\neq l_2=l_4} \cdot  +
  \sum_{l_1=l_4\neq l_2=l_3} \cdot    + \sum_{l_1=l_2=l_3=l_4} \cdot
\eqan
Applying these summations to the three  moments leads to:
\bqan
&&  \expec{\vec{q}}{\sum_{l} q_{l}^2 } \expec{X} {\left( 1_{X=k_1} - \theta_{k_1} \right)\left( 1_{X=k_2} -  \theta_{k_2} \right) } \\
&&  \expec{\vec{q}}{\sum_{l} q_{l}^3 } \expec{X} { \left( 1_{X=k_1} - \theta_{k_1} \right)\left(  1_{X=k_2} - \theta_{k_2} \right) \left( 1_{X=k_3} -  \theta_{k_3} \right)} \\
 &&  \expec{\vec{q}}{\sum_{l} q_{l}^4} \expec{X} {\left( 1_{X=k_1}- \theta_{k_1} \right)\left(  1_{X=k_2} - \theta_{k_2} \right) \left( 1_{X=k_3} -  \theta_{k_3}  \right)
         \left( 1_{X=k_4} - \theta_{k_4} \right)} \\
    && ~~   +  \left( \left( \expec{\vec{q}}{\sum_{l} q_{l}^2} \right)^2 -   \expec{\vec{q}}{\sum_{l} q_{l}^4} \right)   \\
  && ~~~~~ \Big( \expec{X} {\left(   1_{X=k_1} - \theta_{k_1} \right)\left(  1_{X=k_2}  - \theta_{k_2} \right)  }
     \expec{X} {\left(   1_{X=k_3} - \theta_{k_3} \right)\left(  1_{X=k_4} - \theta_{k_4}  \right) }\\
    && ~~  ~~~  ~~
    \expec{X} {\left(  1_{X=k_1} - \theta_{k_1} \right)\left(  1_{X=k_3} - \theta_{k_3}  \right)  }
     \expec{X} {\left( 1_{X=k_2}  -  \theta_{k_2} \right)\left( 1_{X=k_4} - \theta_{k_4} \right) }\\
    && ~~~~~  ~~
    \expec{X} {\left( 1_{X=k_1} - \theta_{k_1} \right)\left( 1_{X=k_4} - \theta_{k_4}  \right)  }
     \expec{X} {\left(  1_{X=k_2} -  \theta_{k_2} \right)\left(  1_{X=k_3} - \theta_{k_3} \right) } \Big)
\eqan
The expected sum of powers of $\vec{q}$ we solve for below. The
expectation of $X$ is for the multivariate discrete (or a
multinomial with $N=1$), so the values are known for the various
cases of $k_1,k_2,...$. For example, when $k_1 \neq k_2$, $
\expec{X} {\left( 1_{X=k_1} - \theta_{k_1}\right)\left( 1_{X=k_2}
- \theta_{k_2} \right) }  =-\theta_{k_1}\theta_{k_2}$.

The expected sum of powers of $\vec{q}$ is obtained as follows. From
first principles, it can be seen that
\beqn
  \expec{\vec{q}}{M}\,=\,\expec{\vec{q}}{1-(1-q_k)^N}
\eeqn
For $N=2$ and rearranging terms we get
\beqn
  \expec{\vec{q}}{M_2}~=~2-\expec{\vec{q}}{\sum_l q_l^2}
\eeqn
Applying Lemma~\ref{lem-part} one gets the closed form expression
for the left-hand side. Likewise, we get:
\bqan
  \expec{\vec{q}}{\sum_l q_l^2} &=& \frac{1-a}{1+b} \\
  \expec{\vec{q}}{\sum_l q_l^3} &=& \frac{(1-a)(2-a)}{(1+b)(2+b)} \\
  \expec{\vec{q}}{\sum_l q_l^4} &=& \frac{(1-a)(2-a)(3-a)}{(1+b)(2+b)(3+b)} \\
  \expec{\vec{q}}{\sum_l q_l^5} &=& \frac{(1-a)(2-a)(3-a)(4-a)}{(1+b)(2+b)(3+b)(4+b)}
\eqan
Combining the resultant formula yields the cases in the lemma.
\qed

\end{document}